\documentclass[a4paper,11pt]{article}
\usepackage{tipa}
\usepackage{bbm}
\usepackage{bbding}

\usepackage{stmaryrd}
\usepackage{amsfonts}
\usepackage{graphicx}
\usepackage{epsfig}
\usepackage{amssymb}
\usepackage{amssymb,amsmath,amsthm}
\usepackage{pifont}
\usepackage{graphicx}
\usepackage{amsmath}
\usepackage{amssymb,amsmath,amsthm}
\usepackage{indentfirst}
\large\normalsize
\def\dis{\displaystyle}

\numberwithin{equation} {section} \makeatletter
\renewcommand{\@seccntformat}[1]{\csname
the#1\endcsname.\hspace{0.5em}} \makeatother

\begin{document}

 \title{Flux-approximation limits of solutions to the Brio system with two independent parameters
 \thanks{\ Supported by National Natural Science Foundation of China (11501488), the Scientific Research
 Foundation of Xinyang Normal University(No. 0201318), Nan Hu Young Scholar Supporting Program of XYNU,
 Yunnan Applied Basic Research Projects (2018FD015) and the Scientific Research Foundation Project
of Yunnan Education Department (2018JS150).}
 }
\author{ Yanyan Zhang $^{1}$\thanks{\ Corresponding author. E-mail address: zyy@xynu.edu.cn (Y. Zhang), yuzhang13120@126.com (Y. Zhang).}
  \ \ Yu Zhang $^{2}$ \\
{\footnotesize\slshape{$^{1}$College of Mathematics and Statistics,  Xinyang Normal
University, Xinyang 464000, PR China}}\\
{\footnotesize\slshape{$^{2}$Department of Mathematics,  Yunnan Normal
University, Kunming 650500, PR China}}
}

\date{}
\maketitle

\textbf{Abstract:} By the flux-approximation method, we
study limits of Riemann solutions to the Brio system with two independent parameters.
The Riemann problem of the perturbed system is solved analytically,
and four kinds of solutions are obtained constructively.
It is shown that, as the two-parameter flux perturbation vanishes, any two-shock-wave
and two-rarefaction-wave solutions of the perturbed Brio system converge to
the delta-shock and vacuum solutions of the transport equations, respectively.
In addition, we specially pay attention to the Riemann problem of a simplified system of conservation laws
derived from the perturbed Brio system by neglecting some quadratic term.
As one of the purebred parameters of the Brio system goes to zero, the
solution of which consisting of two shock waves tends to a delta-shock solution to this simplified system.
By contrast, the solution containing two rarefaction waves converges to a
contact discontinuity and a rarefaction wave
of the simplified system. What is more, the formation mechanisms of delta shock waves
under flux approximation with both two parameters and only one parameter are clarified.

\textbf{Keywords:}  Brio system; Transport equations; Riemann problem; Delta shock wave; Vacuum;
Flux approximation.

\date{}

\section{Introduction }

As to our knowledge, in the past over two decades, the
delta shock wave has been systematically studied
by a large number of scholars. For example, see the results
in \cite{Korchinski, Keyfitz, Huang, Shelkovich, Tan1, Tan2, Sheng, Yang} and the references cited therein. Particularly,
in the related researches of delta shock waves,
one very interesting topics is to explore the formation of
delta shock waves and vacuum states in solutions,
which correspond to the phenomena of concentration and cavitation, respectively.
At this moment, an effective approach is to use the so called vanishing pressure limit method,
which was early proposed by Chen and Liu \cite{CHL1, CHL2}
to study the formation of delta shock waves and vacuums for the Euler equations of
isentropic and nonisentropic gas dynamics, respectively. See also Li \cite{Li}
for the isothermal Euler equations with zero temperature.

In view that the vanishing pressure limit method is only focused on pressure perturbation, for more general of
physical consideration, Yang and Liu
recently  introduced the flux-approximation method \cite{yang-liu1} to study the limit behavior
of solutions to the isentropic Euler equations with flux perturbation.
The flux-approximation method, generally speaking,  is a natural generalization
of the vanishing pressure limit method.
The main idea of it is to introduce some
small perturbed parameters in the flux function of the system, then
the limits of solutions to the perturbed system can be studied
by taking the perturbed parameters go to zero. This method
has been successfully applied to study the formation of delta shock waves, say, Yang and Liu
\cite{yang-liu2}  for the nonisentropic fluid flows, Yang and Zhang \cite{yang-zhangyu1, yang-zhangyu2}
for the relativistic Euler equations, Sun \cite{Sun1} for the transport equations, etc.

Motivated by the works mentioned above, we in this paper introduce the
following system of conservation laws
\begin{align}\label{eq1.1}
\left\{\begin{array}{l}
u_{t}+(\frac{1}{2}u^{2}+\frac{1}{2}\epsilon_{1}v^{2})_{x}=0, \cr\noalign {\vskip2truemm}
v_{t}+(uv-\epsilon_{2}v)_{x}=0,
\end{array}\right.
\end{align}
where $\epsilon_{1}$, $\epsilon_{2}>0$ are two independent parameters.
Here, we are pleased to mention that, even though the parameters $\epsilon_{1}$ and $\epsilon_{2}$ are considered
very small and govern the strength of flux, they do not vanish in general.
We just propose to include these small parameters
to explore the limit behaviors of solutions to the system \eqref{eq1.1}.
Introducing some small perturbed parameters in the
flux function of the system, on one hand, it
can be used to control some dynamical behaviors of fluids \cite{yang-liu1} from a physical point of view;
on the other hand, it can be used to numerical calculation \cite{Smith, RJ, Cap, Colombeau, Dedovic} because it can perturb
the non-strictly hyperbolic system into a nearby strictly hyperbolic one from a mathematical point of view.

System \eqref{eq1.1} is called the perturbed Brio system, since it can be viewed as
a perturbed model of the Brio system \cite{Brio}
\begin{align}\label{eq1.22}
\left\{\begin{array}{l}
u_{t}+\left(\frac{u^{2}+v^{2}}{2}\right)_{x}=0, \cr\noalign {\vskip2truemm}
v_{t}+(uv-v)_{x}=0,
\end{array}\right.
\end{align}
which arises as a simplified model in ideal magnetohydrodynamics (MHD)
and corresponds to coupling the fluid dynamic equations with Maxwell's equations of electrodynamics.
The Brio system \eqref{eq1.22} is strictly hyperbolic and genuinely nonlinear at $\{(u, v): u\in R, v>0\}$ and
$\{(u, v): u\in R, v<0\}$, but not on the whole of $R^{2}$.
In \cite{Hayes}, Hayes et al found that the genuine nonlinearity was lost when $v$ changed sign.
For solutions crossing the line $v=0$, delta-shock solution might have to be used,
and non-uniqueness of Riemann solutions was anticipated.

Letting $\epsilon_{1}, \epsilon_{2}\rightarrow 0$,
the system \eqref{eq1.1} transforms formally into the following transport equations
\begin{align}\label{eq1.2}
\left\{\begin{array}{l}
u_{t}+(\frac{1}{2}u^{2})_{x}=0, \cr\noalign {\vskip2truemm}
v_{t}+(uv)_{x}=0,
\end{array}\right.
\end{align}
where $u$ is velocity and $v$ the density.
The system \eqref{eq1.2} has a physical context
and describes some important physical
phenomena. It can be used to  model the motion of free particles which stick under
collision \cite{Brenier}, and to describe the formation of large-scale structures of the
universe \cite{Shandarin, Ya}. See also \cite{Berthelin, Edwards} for more related applications.
It is non-strictly hyperbolic, with a linearly degenerate characteristic field, and
has been studied in the numerous papers such as \cite{Joseph, Shen2, Bouchut, Sheng, Huang}.
Interestingly, the delta shock wave and vacuum state do appear in the Riemann solutions.

However, if we only take  $\epsilon_{1}\rightarrow 0$, then \eqref{eq1.1}
becomes the following single-parameter-perturbation model
\begin{align}\label{eq1.5}
\left\{\begin{array}{l}
u_{t}+(\frac{1}{2}u^{2})_{x}=0, \cr\noalign {\vskip2truemm}
v_{t}+(uv-\epsilon_{2}v)_{x}=0,
\end{array}\right.
\end{align}
which was used to approximate the transport equation by Shen etc \cite{Shen1}.
They proved that, as $\epsilon_{2} \rightarrow 0$, the Riemann solutions of \eqref{eq1.5}
converge to those of the transport equations \eqref{eq1.2}.
Specifically, as $\epsilon_2=1$, system \eqref{eq1.5} is called the simplified
Brio system, which was derived in \cite{Hayes} by neglecting the
$v^{2}$ term in the flux function of the first equation of Brio system \eqref{eq1.22}.
What is more, it is found that the delta shock wave occurs in the solutions of \eqref{eq1.5}.

The main purpose of this paper is to discuss the formation of delta shock waves and
vacuum states in the vanishing flux-approximation limit of solutions to the
 perturbed Brio model \eqref{eq1.1}. It remains to be seen whether
or not the limits $\epsilon_{1}, \epsilon_{2} \rightarrow 0$ and $\epsilon_{1} \rightarrow 0$
of solutions to the Riemann problem for the perturbed Brio system \eqref{eq1.1} are identical with those
for the transport equations \eqref{eq1.2} and the single-parameter-perturbation model \eqref{eq1.5}, respectively.
In what follows, we outlook the context of each section of this paper.

In Section 2,  we solve  delta shock waves and vacuum states for transport equations  \eqref{eq1.2}.

Section 3 deals with the Riemann problem for the single-parameter-perturbation system \eqref{eq1.5} with initial data
\begin{align}\label{eq1.6}
(u,v)(0,x)=\  \left\{\begin{array}{lll}
 (u_{-},v_{-}), & x <0,\cr\noalign{\vskip2truemm}
(u_{+},v_{+}), &  x>0,
\end{array}\right.
\end{align}
where $u_{\pm}$ and $v_{\pm}$ are arbitrary constants. Two kinds of Riemann solutions consisting
of rarefaction waves, shock waves and contact discontinuities
are constructed when $u_{-}-u_{+}<2\epsilon_{2}$.
While when $u_{-}-u_{+}>2\epsilon_{2}$, the delta shock wave appears in solutions.
Nevertheless, the velocity of delta shock wave in system \eqref{eq1.5} is quite different from that of
\eqref{eq1.2}, since it will no longer equal to the value of $u$ on the discontinuity line,
which implies that the flux perturbation works in the transport equations \eqref{eq1.2}.

In Section 4, we investigate the Riemann problem of the perturbed Brio system \eqref{eq1.1}.
For convenience, we only consider the case when $v>0$,
that is, the system \eqref{eq1.1} is strictly hyperbolic and genuinely nonlinear,
which also means that the perturbation adopted here transforms the non-strictly hyperbolic
system \eqref{eq1.2} into the strictly hyperbolic one.  By analyzing in phase plane,
we construct four kinds of Riemann solutions with the
classical waves involving rarefaction waves and shock waves.

 Sections 5 studys the flux-approximation limits of
solutions to \eqref{eq1.1}, \eqref{eq1.6} as both parameters go to zero.
Concretely, it is rigorously proved that, as  $\epsilon_{1}, \epsilon_{2}\rightarrow 0$,
any two-shock Riemann solution to the perturbed  Brio system \eqref{eq1.1}
tends to a delta-shock solution to the transport equations \eqref{eq1.2},
and the intermediate state between the two shock waves tends to a weighted $\delta$-measure that forms a delta shock wave.
Meanwhile, it is also shown that, any two-rarefaction-wave Riemann solution to the perturbed  Brio system \eqref{eq1.1}
converges to a two-contact-discontinuity solution to the transport equations \eqref{eq1.2},
whose non-vacuum intermediate state in between tends to a vacuum state.
These results present that the delta-shock and vacuum solutions of \eqref{eq1.2}
can be obtained by a flux-approximation limit of solutions to the perturbed  Brio system \eqref{eq1.1}.

In Section 6, we discuss the limit behaviors of Riemann solutions to the  perturbed  Brio system \eqref{eq1.1}
as $\epsilon_{1}\rightarrow 0$. It is shown that the Riemann solutions of \eqref{eq1.1} converge
to the corresponding Riemann solutions of the
single-parameter-perturbation model \eqref{eq1.5}. Especially,
any two-shock-wave Riemann solution of \eqref{eq1.1}
tends to a delta-shock solution to the perturbed model \eqref{eq1.5}.

Following the above analysis, from the point of hyperbolic conservation laws,
the above two different convergence processes show the
two kinds of occurrence mechanism on the formation of
delta shock wave. When  $\epsilon_{1}$ and $\epsilon_{2}$ decrease simultaneous,
the strict hyperbolicity of the limiting system fails(see Section 5), which leads to the formation
of delta shock wave. This point is the same as that in \cite{Shen1, Sun1}.
While when only one parameter $\epsilon_{1}$ decreases,
the strict hyperbolicity of
the limiting system is preserved(see Section 6), but the delta shock wave
still occurs, which is different from that in \cite{Shen1, Sun1}.
In addition, the above results also indicate the fact that
different flux approximations have their respective effects on the formation of
delta shock waves.

\section{Preliminaries}

In this section, let us recall the Riemann solutions to the simplified transport equations  \eqref{eq1.2}.
More details can be found in \cite{Joseph, Shen2, Sheng}.

The system  \eqref{eq1.2} has duplicate eigenvalues $\lambda=u$ with corresponding right eigenvectors
$r=(1,0)^{T}$. Since  $\nabla\lambda \cdot \vec{r}=0$, so \eqref{eq1.2}
is full linear degenerate and elementary waves involve only contact discontinuities.
The left state $(u_{-},v_{-})$ and right state $(u_{+},v_{+})$
can be connected by classical elementary waves(contact discontinuities and vacuum) or a delta shock wave.
Depending on the choice of initial data, there are two possible wave patterns for
solutions of Riemann problem \eqref{eq1.2} and \eqref{eq1.6}.

When $u_{-}<u_{+}$, the Riemann solution consists of two contact discontinuities $J$ with a vacuum in between,
which can be shown as
\begin{align}
(u, v)(t,x)=\left\{\begin{array}{lll}
 (u_{-},v_{-}), & -\infty<x <u_{-}t ,\cr\noalign{\vskip2truemm}
 (\frac{x}{t}, 0), & u_{-}t\leqslant x\leqslant u_{+}t,\cr\noalign{\vskip2truemm}
(u_{+},v_{+}), &  u_{+}t <x<+\infty.
\end{array}\right.
\end{align}

When $u_{-}>u_{+}$, the Riemann solution cannot be constructed by using the classical waves,
and the delta shock wave appears. To define this kind of solution, the following weighted delta function supported on a curve should
be introduced.

A two-dimensional weighted delta function
$w(s)\delta_S$ supported on a smooth curve $S$ parameterized as
$t=t(s), x=x(s)(a\leqslant s\leqslant b)$ can be defined by
\begin{align}\label{eq2.2}
\Big\langle w(t(s))\delta_S,
\varphi(t,x)\Big\rangle=\int_a^bw(t(s))\varphi\big(t(s),x(s)\big)\mbox{d}s
\end{align}
for all test functions $\varphi\in C^\infty_0([0,+\infty)\times(-\infty, +\infty))$.

Based on this definition, a delta-shock solution of  \eqref{eq1.2}  can be represented in the
following form
\begin{align}\label{eq2.3}
u(t,x)=u_{0}(x,t),\ \ \ \ \ \ \  v(t,x)=v_{0}(t,x)+w(t)\delta_{s},
\end{align}
where $S=\{(t,\sigma t):0\leqslant t <\infty\}$, and
\begin{align}
\begin{array}{l}
u_{0}(t,x)=u_{-}+[u]H(x-\sigma t), \ \  v_{0}(t,x)=v_{-}+[v]H(x-\sigma t),  \ \   w(t)=\sigma[v]-[uv],
\end{array}
\end{align}
in which $[G]=G_{+}-G_{-}$ expresses the jump of the quality $G$ across the curve $S$,
$\sigma$ is the velocity of the delta shock wave, and $H(x)$  the Heaviside function.

As mentioned in \cite{Sheng} that the solution $(u,v)$ constructed above satisfies that
\begin{align}
\begin{array}{l}
\langle  u, \varphi_{t}\rangle+\langle \frac{1}{2}u^{2}, \  \varphi_{x}\rangle=0,\cr\noalign{\vskip2truemm}
\langle v, \varphi_{t}\rangle+\langle uv, \ \varphi_{x}\rangle=0,
\end{array}
\end{align}
for all test functions $\varphi\in C^\infty_{0}\big([0,+\infty)\times(-\infty,+\infty)\big)$, where
$$\langle v, \varphi \rangle=\dis\int_{0}^{+\infty}\dis\int_{-\infty}^{+\infty}v_{0}\varphi \mbox{d}x\mbox{d}t
+\langle w\delta_{S}, \varphi \rangle,$$
$$\langle  uv, \varphi \rangle=\dis\int_{0}^{+\infty}\dis\int_{-\infty}^{+\infty}u_{0}v_{0}\varphi \mbox{d}x\mbox{d}t
+\langle \sigma w\delta_{S}, \varphi \rangle.$$

Then, a unique solution of \eqref{eq1.2} and \eqref{eq1.6}
 containing a weighted $\delta$-measure
is given as
\begin{align}\label{eq2.6}
(u, v)(t,x)=\left\{\begin{array}{lll}
 (u_{-},v_{-}), &  x <x(t),         \cr\noalign{\vskip2truemm}
 (\sigma, w(t)\delta(x-\sigma t)), & x=x(t),              \cr\noalign{\vskip2truemm}
(u_{+},v_{+}),  & x>x(t),
\end{array}\right.
\end{align}
in which $x(t)$, $\sigma$ and  $w(t)$ satisfy the generalized Rankine-Hugoniot relation
\begin{align}\label{eq2.7}
\left\{\begin{array}{l}
\dfrac{\mbox{d}x}{\mbox{d}t}=\sigma,\cr\noalign{\vskip2truemm}
\dfrac{ \mbox{d}w(t)}{\mbox{d}t}=[v]\sigma-[uv],\cr\noalign{\vskip2truemm}
[u]\sigma=\Big[\dfrac{u^{2}}{2}\big].
\end{array}\right.
\end{align}

In addition, the entropy condition
$$
u_{+}<\sigma<u_{-}
$$
should be supplemented in order to guarantee the uniqueness.

By solving the generalized Rankine-Hugoniot relation \eqref{eq2.7} with the initial data
$w(0)=0$, $x(0)=0$, one can obtain that
\begin{align}\label{eq2.8}
\sigma=\frac{1}{2}(u_{-}+u_{+}),  \ \ \ w(t)=\frac{1}{2}(v_{-}+v_{+})(u_{-}-u_{+})t.
\end{align}
Thus, the delta-shock solution defined by \eqref{eq2.6} with \eqref{eq2.8} is obtained.

\section{Riemann solutions of system \eqref{eq1.5}}
In this section, we shall solve the Riemann problem of \eqref{eq1.5} and \eqref{eq1.6}.

System \eqref{eq1.5} has two eigenvalues
$$
\begin{array}{l}
\lambda^{\epsilon_{2}}_{1}=u-\epsilon_{2}, \hspace{2cm}
\lambda^{\epsilon_{2}}_{2}=u
\end{array}
$$
with the corresponding right eigenvectors
$$
\begin{array}{l}
\overrightarrow{r_{1}}^{\epsilon_{2}}=(0, 1)^{T},\hspace{1.5cm}
 \overrightarrow{r_{2}}^{\epsilon_{2}}=(\epsilon_{2}, v)^{T}.
\end{array}
$$
Therefore, the  system  is strictly hyperbolic for
$\epsilon_{2}>0$. In addition,
it is easily to check that
$\nabla \lambda^{\epsilon_{2}}_{1}\cdot \overrightarrow{r_{1}}^{\epsilon_{2}}=0$
and $\nabla \lambda^{\epsilon_{2}}_{2}\cdot \overrightarrow{r_{2}}^{\epsilon_{2}}=\epsilon_2$,
which mean that $\lambda^{\epsilon_{2}}_{1}$
is always linearly degenerate and $\lambda^{\epsilon_{2}}_{2}$ is genuinely nonlinear for
$\epsilon_{2}>0$.

Seeking the self-similar solution
$$
\begin{array}{l}
(u,v)(t,x)=(u,v)(\xi), \ \ \ \xi=x/t.
\end{array}
$$
Equivalently, the Riemann problem \eqref{eq1.5} and \eqref{eq1.6} is reduced to
\begin{align}
\left\{\begin{array}{l}
-\xi u_{\xi}+(\frac{1}{2}u^{2})_{\xi}=0, \cr\noalign {\vskip2truemm}
-\xi v_{\xi}+(uv-\epsilon_{2}v)_{\xi}=0 \cr\noalign {\vskip2truemm}
(u, v)(0,\pm\infty)=(u_{\pm},v_{\pm}),
\end{array}\right.
\end{align}
which provides either the constant state solution or the rarefaction wave
\begin{align}
R(u_{-},v_{-}):\left\{\begin{array}{l}
 \xi =\lambda^{\epsilon_{2}}_{2}=u, \cr\noalign {\vskip2truemm}
\epsilon_{2}\mbox{ln}v-u=\epsilon_{2}\mbox{ln}v_{-}-u_{-}, \cr\noalign {\vskip2truemm}
u>u_{-}.
\end{array}\right.
\end{align}
For a bounded discontinuity at $\xi=w$, the following Rankine-Hugoniot relation
\begin{align}
\left\{\begin{array}{l}
 w[u]=\big[\frac{1}{2}u^{2}], \\[0.5cm]
w[v]=\big[uv-\epsilon_{2}v]
\end{array}\right.
\end{align}
holds. When $[u]\neq0$,  we get the shock wave
\begin{align}
S(u_{-},v_{-}):\left\{\begin{array}{l}
 w_{1} =\dfrac{u_{-}+u}{2}, \cr\noalign {\vskip2truemm}
\dfrac{v}{v_{-}}=\dfrac{u-u_{-}+2\epsilon_{2}}{u_{-}-u+2\epsilon_{2}},\cr\noalign {\vskip2truemm}
u<u_{-}<u+2\epsilon_{2}.
\end{array}\right.
\end{align}
While when $[u]=0$, it corresponds to a contact discontinuity
\begin{align}
J: \ \ \ \ w^{\epsilon_{2}}=u_{-}-\epsilon_{2}=u_{+}-\epsilon_{2}.
\end{align}

However, when $u_{-}-u_{+}>2\epsilon_{2}$, the solution cannot be constructed by classical waves.
At this moment, the delta-shock solution containing Dirac delta function in the state variables $v$
will be considered.
Under the definition \eqref{eq2.2}, we give the definition of the delta-shock solution to \eqref{eq1.5} and \eqref{eq1.6}.

\vspace{0.3cm}

\noindent{ \bf Definition 3.1} \ A pair $(u,v)$ is called a delta shock wave type solution of \eqref{eq1.5} in the sense of distributions
if there exist a smooth curve $S$ and the function $w^{\epsilon_{2}}(t)$ such that $u$ and $v$ are represented in the
following form
\begin{align}
 u=\bar{u}(t,x), \ \ \ \ \ \ v=\bar{v}(t,x)+w^{\epsilon_{2}}(t)\delta_{S},
\end{align}
where $\bar{u}, \bar{v} \in L^{\infty}(R\times (0,+\infty); R)$, $w^{\epsilon_{2}}(t) \in C^{1}(S)$, $u|_{S}=u^{\epsilon_{2}}_{\delta}=\sigma^{\epsilon_{2}}+\epsilon_{2}$,
$\sigma^{\epsilon_{2}}$ is the tangential derivative of curve $S$,  and they satisfy
\begin{align}\label{eq3.7}
\begin{array}{l}
\langle  u, \varphi_{t}\rangle+\langle \frac{1}{2}u^{2}, \  \varphi_{x}\rangle=0, \cr\noalign{\vskip2truemm}
\langle v, \varphi_{t}\rangle+\langle uv-\epsilon_{2}v, \ \varphi_{x}\rangle=0
\end{array}
\end{align}
for all test functions $\varphi\in C^\infty_{0}\big((-\infty,+\infty)\times [0,+\infty)\big)$, where
$$\langle v, \varphi \rangle=\dis\int_{0}^{+\infty}\dis\int_{-\infty}^{+\infty}\bar{v}\varphi \mbox{d}x\mbox{d}t
+\langle w^{\epsilon_{2}}(t)\delta_{S}, \varphi \rangle,$$
$$\langle uv, \varphi \rangle=\dis\int_{0}^{+\infty}\dis\int_{-\infty}^{+\infty}\bar{u}\bar{v}\varphi \mbox{d}x\mbox{d}t
+\langle w^{\epsilon_{2}}(t)u^{\epsilon_{2}}_{\delta}\delta_{S}, \varphi \rangle.$$

Using this definition, we can define the delta-shock solution
of \eqref{eq1.5}  with the discontinuity $x=x(t)$ of the form
\begin{align}\label{eq3.8}
(u,v)(t,x)=\  \left\{\begin{array}{lll}
 (u_{-},v_{-}), &  x<x(t),\cr\noalign{\vskip2truemm}
(u^{\epsilon_{2}}_{\delta}, w^{\epsilon_{2}}(t)\delta(x-x(t))), &   x=x(t),\cr\noalign{\vskip2truemm}
(u_{+},v_{+}), &  x>x(t),
\end{array}\right.
\end{align}
where $(u_{-},v_{-})$ and $(u_{+},v_{+})$ are piecewise smooth solutions of
\eqref{eq1.5}, $\delta(\cdot)$ is the standard Dirac measure, $x(t)\in C^{1}$,
 $w^{\epsilon_{2}}(t)$ is strength of delta shock wave
 and $u^{\epsilon_{2}}_{\delta}$ is the corresponding value of $u$
on the line $x=x(t)$.

The solution $(u,v)(t,x)$ defined in \eqref{eq3.8} satisfies \eqref{eq3.7} in the sense of distributions if
it satisfies the relation
\begin{align}\label{eq3.9}
\left\{\begin{array}{l}
\dfrac{\mbox{d}x}{\mbox{d}t}=\sigma^{\epsilon_{2}},\cr\noalign{\vskip2truemm}
\dfrac{ \mbox{d}w^{\epsilon_{2}}(t)}{\mbox{d}t}=[v]\sigma^{\epsilon_{2}}-[uv-\epsilon_{2}v],\cr\noalign{\vskip2truemm}
[u]\sigma^{\epsilon_{2}}=\Big[\dfrac{u^{2}}{2}\Big]
\end{array}\right.
\end{align}
and
\begin{align}\label{eq3.10}
u|_{x=x(t)}=u^{\epsilon_{2}}_{\delta}=\sigma^{\epsilon_{2}}+\epsilon_{2}.
\end{align}

In fact, for any test functions $\varphi\in C^\infty_{0}\big((-\infty,+\infty)\times [0,+\infty)\big)$,
if  \eqref{eq3.9} and \eqref{eq3.10} hold, then by Green's formulation and integrating by parts, it yields that
$$
\begin{array}{l}
\langle v, \ \varphi_{t}\rangle+\langle uv-\epsilon_{2}v, \ \varphi_{x}\rangle\cr\noalign{\vskip4truemm}
 =\dis \int_{0}^{+\infty}\int_{-\infty}^{x(t)}\big(v_{-}\varphi_{t}+v_{-}(u_{-}-\epsilon_{2})\varphi_{x}\big)dxdt
 +\dis \int_{0}^{+\infty}\int_{x(t)}^{+\infty}\big(v_{+}\varphi_{t}+v_{+}(u_{+}-\epsilon_{2})\varphi_{x}\big)dxdt
 \cr\noalign{\vskip4truemm}
\hspace{0.5cm} + \dis \int_{0}^{+\infty}w^{\epsilon_{2}}(t)(\varphi_{t}+(u^{\epsilon_{2}}_{\delta}-\epsilon_{2})\varphi_{x})dt
 \cr\noalign{\vskip4truemm}
=-\dis\oint-(v_{-}(u_{-}-\epsilon_{2})\varphi)dt+v_{-}\varphi dx+\dis\oint-(v_{+}(u_{+}-\epsilon_{2})\varphi)dt+v_{+}\varphi dx
 \cr\noalign{\vskip4truemm}
\hspace{0.5cm}+\dis \int_{0}^{+\infty}w^{\epsilon_{2}}(t)(\varphi_{t}+(u^{\epsilon_{2}}_{\delta}-\epsilon_{2})\varphi_{x})dt \cr\noalign{\vskip4truemm}
=\dis \int_{0}^{+\infty}(v_{-}(u_{-}-\epsilon_{2})-v_{+}(u_{+}-\epsilon_{2}))\varphi dt+\dis \int_{0}^{+\infty}(v_{+}-v_{-})\sigma^{\epsilon_{2}}\varphi dt-\dis \int_{0}^{+\infty}\frac{dw^{\epsilon_{2}}(t)}{dt}\varphi dt \cr\noalign{\vskip4truemm}
=\dis \int_{0}^{+\infty}\Big(\sigma^{\epsilon_{2}}[v]-[uv-\epsilon_{2}v]-\frac{dw^{\epsilon_{2}}(t)}{dt}\Big)\varphi dt\cr\noalign{\vskip4truemm}
=0.
\end{array}
$$
That is, the second equation of \eqref{eq3.7} hold. The rest one can be checked  in a similar way.

Equations \eqref{eq3.9} and \eqref{eq3.10} are called the generalized Rankine-Hugoniot relation of delta shock waves.
Furthermore, the entropy condition
\begin{align}
u_{-}-\epsilon_{2}>\sigma^{\epsilon_{2}}>u_{+}+\epsilon_{2}
\end{align}
should be supplemented to guarantee the uniqueness.

By solving the ordinary differential
equations \eqref{eq3.9} and \eqref{eq3.10} with initial data
$t=0: \ x(0)=0, \ \  w^{\epsilon_{2}}(0)=0$,
we get that
\begin{align}\label{eq3.12}
\begin{array}{l}
\sigma^{\epsilon_{2}}=u^{\epsilon_{2}}_{\delta}-\epsilon_{2}=\frac{1}{2}(u_{-}+u_{+}),  \cr\noalign {\vskip2truemm} w^{\epsilon_{2}}(t)=\frac{1}{2}(v_{+}(u_{-}-u_{+}+2\epsilon_{2})-v_{-}(u_{+}-u_{-}+2\epsilon_{2}))t.
\end{array}
\end{align}
Given a constant state $(u_{-},v_{-})$, these four wave curves divide the half-phase plane into three regions
\begin{align}
\begin{array}{l}
I=\{(u,v)\mid u_{-}<u<+\infty\}; \ \ \ II=\{(u,v)\mid u_{-}-2\epsilon_{2}<u<u_{-}\};\cr\noalign{\vskip4truemm}
III=\{(u,v)\mid -\infty<u<u_{-}-2\epsilon_{2}\}.
\end{array}
\end{align}

According to the state $(u_{+},v_{+})$ in the different regions, the solution is $J+R$ when $(u_{+},v_{+})\in I(u_{-},v_{-}) $,
$J+S$ when $(u_{+},v_{+})\in II(u_{-},v_{-}) $, and delta shock wave when $(u_{+},v_{+})\in III(u_{-},v_{-}) $(see Fig. 1),
where the symbol $``+" $ means ``followed by".

\section{Riemann solutions of the perturbed Brio system \eqref{eq1.1}}

Now, we construct the Riemann solutions to \eqref{eq1.1} and \eqref{eq1.6}.
The eigenvalues of \eqref{eq1.1} are
$$
\begin{array}{l}
\lambda^{\epsilon_{1}\epsilon_{2}}_{1}=u-\dfrac{1}{2}\epsilon_{2}-\dfrac{\sqrt{\epsilon^{2}_{2}+4\epsilon_{1}v^{2}}}{2}, \ \ \
\lambda^{\epsilon_{1}\epsilon_{2}}_{2}=u-\dfrac{1}{2}\epsilon_{2}+\dfrac{\sqrt{\epsilon^{2}_{2}+4\epsilon_{1}v^{2}}}{2}
\end{array}
$$
with $\lambda^{\epsilon_{1}\epsilon_{2}}_{1}<\lambda^{\epsilon_{1}\epsilon_{2}}_{2}$,
so the system  \eqref{eq1.1} is strictly hyperbolic for $\epsilon_{1}, \epsilon_{2}>0$. The corresponding right eigenvectors are
given by
$$
\begin{array}{l}
\overrightarrow{r_{1}}^{{\epsilon_{1}\epsilon_{2}}}=\Big(\dfrac{1}{2}\epsilon_{2}-\dfrac{\sqrt{\epsilon^{2}_{2}+4\epsilon_{1}v^{2}}}{2},\ \ v \Big)^{T},  \ \ \
 \overrightarrow{r_{2}}^{{\epsilon_{1}\epsilon_{2}}}=\Big(\dfrac{1}{2}\epsilon_{2}+\dfrac{\sqrt{\epsilon^{2}_{2}+4\epsilon_{1}v^{2}}}{2}, \ \ v \Big)^{T}.
\end{array}
$$
Checking genuine nonlinearity for the perturbed Brio system \eqref{eq1.1},
 we find that
\begin{align}\label{eq4.1}
\begin{array}{l}
\nabla \lambda^{\epsilon_{1}\epsilon_{2}}_{i}\cdot \overrightarrow{r_{i}}^{{\epsilon_{1}\epsilon_{2}}}=\epsilon_{1}v\Big(2\pm\dfrac{\epsilon_{2}}
{\sqrt{\epsilon^{2}_{2}+4\epsilon_{1}v^{2}}}\Big), \ \ i=1,2.
\end{array}
\end{align}
Therefore, both two characteristic fields are genuinely nonlinear when $v>0$ and $\epsilon_{1}, \epsilon_{2}>0$.

As before, we look for the self-similar solution, then \eqref{eq1.1} becomes
\begin{align}\label{eq4.2}
\left\{\begin{array}{l}
-\xi u_{\xi}+(\frac{1}{2}u^{2}+\frac{1}{2}\epsilon_{1}v^{2})_{\xi}=0, \cr\noalign {\vskip2truemm}
-\xi v_{\xi}+(uv-\epsilon_{2}v)_{\xi}=0
\end{array}\right.
\end{align}
with the boundary condition
\begin{align}
\begin{array}{l}
(u, v)(0,\pm\infty)=(u_{\pm},v_{\pm}).
\end{array}
\end{align}

For the smooth solution, we write \eqref{eq4.2} in matrix form as
$$
\left(\begin{array}{ccc}
 -\xi+u& \epsilon_{1}v \\
v \ \ & -\xi+u-\epsilon_{2}
\end{array}\right)
\left(\begin{array}{ccc}
 u \\ v
\end{array}\right)_{\xi}=0.
$$
Besides the constant state
$$
(u, v)(\xi)=\mbox{constant},
$$
it provides the backward rarefaction wave
\begin{align}\label{eq4.4}
R_{1}(u_{-},v_{-}):\left\{\begin{array}{l}
 \xi =\lambda^{\epsilon_{1}\epsilon_{2}}_{1}=u-\dfrac{1}{2}\epsilon_{2}-\dfrac{\sqrt{\epsilon^{2}_{2}+4\epsilon_{1}v^{2}}}{2}, \cr\noalign {\vskip2truemm}
 v\mbox{d}u=\Big(\dfrac{1}{2}\epsilon_{2}-\dfrac{\sqrt{\epsilon^{2}_{2}+4\epsilon_{1}v^{2}}}{2}\Big)\mbox{d}v,
\end{array}\right.
\end{align}
or the forward  rarefaction wave
\begin{align}\label{eq4.5}
R_{2}(u_{-},v_{-}):\left\{\begin{array}{l}
 \xi =\lambda^{\epsilon_{1}\epsilon_{2}}_{2}=u-\dfrac{1}{2}\epsilon_{2}+\dfrac{\sqrt{\epsilon^{2}_{2}+4\epsilon_{1}v^{2}}}{2}, \cr\noalign {\vskip2truemm}
 v\mbox{d}u=\Big(\dfrac{1}{2}\epsilon_{2}+\dfrac{\sqrt{\epsilon^{2}_{2}+4\epsilon_{1}v^{2}}}{2}\Big)\mbox{d}v.
\end{array}\right.
\end{align}

From \eqref{eq4.4} and \eqref{eq4.5}, we calculate that
\begin{align}\label{eq4.6}
\begin{array}{l}
 \dfrac{\mbox{d}\lambda_{1}}{\mbox{d}v}=\dfrac{\epsilon_{2}\sqrt{\epsilon^{2}_{2}+4\epsilon_{1}v^{2}}-\epsilon^{2}_{2}-8\epsilon_{1}v^{2}}
 {2v\sqrt{\epsilon^{2}_{2}+4\epsilon_{1}v^{2}}}<0
\end{array}
\end{align}
and
\begin{align}\label{eq4.7}
\begin{array}{l}
 \dfrac{\mbox{d}\lambda_{2}}{\mbox{d}v}=\dfrac{\epsilon_{2}\sqrt{\epsilon^{2}_{2}+4\epsilon_{1}v^{2}}+\epsilon^{2}_{2}+8\epsilon_{1}v^{2}}
 {2v\sqrt{\epsilon^{2}_{2}+4\epsilon_{1}v^{2}}}>0.
\end{array}
\end{align}
Thus, those states which can be connected to $(u_{-},v_{-})$ by $R_{1}$($R_{2}$)
must lie in a direction in which $\lambda^{\epsilon_{1}\epsilon_{2}}_{1}$($\lambda^{\epsilon_{1}\epsilon_{2}}_{2}$) is monotonically
increasing, so the left state $(u_{-},v_{-})$ and right state $(u,v)$
can be connected by $R_{1}$($R_{2}$) with $v<v_{-}$($v>v_{-}$).

Integrating the second equations of \eqref{eq4.4} and \eqref{eq4.5} respectively, one can get that
\begin{align}\label{eq4.8}
R_{1}(u_{-},v_{-}):\left\{\begin{array}{l}
 \xi =\lambda^{\epsilon_{1}\epsilon_{2}}_{1}=u-\dfrac{1}{2}\epsilon_{2}-\dfrac{\sqrt{\epsilon^{2}_{2}+4\epsilon_{1}v^{2}}}{2}, \cr\noalign {\vskip2truemm}
 u-\dfrac{1}{2}\Big(-\sqrt{\epsilon^{2}_{2}+4\epsilon_{1}v^{2}}
+\epsilon_{2}\ln(\sqrt{\epsilon^{2}_{2}+4\epsilon_{1}v^{2}}+\epsilon_{2})\Big)=C_{1},\cr\noalign {\vskip4truemm}
u>u_{-}, \ \ v<v_{-},
\end{array}\right.
\end{align}
and
\begin{align}\label{eq4.9}
R_{2}(u_{-},v_{-}):\left\{\begin{array}{l}
 \xi =\lambda^{\epsilon_{1}\epsilon_{2}}_{2}=u-\dfrac{1}{2}\epsilon_{2}+\dfrac{\sqrt{\epsilon^{2}_{2}+4\epsilon_{1}v^{2}}}{2}, \cr\noalign {\vskip2truemm}
 u-\dfrac{1}{2}\Big(\sqrt{\epsilon^{2}_{2}+4\epsilon_{1}v^{2}}
+\epsilon_{2}\ln(\sqrt{\epsilon^{2}_{2}+4\epsilon_{1}v^{2}}-\epsilon_{2})\Big)=C_{2},\cr\noalign {\vskip4truemm}
u>u_{-}, \ \ v>v_{-},
\end{array}\right.
\end{align}
where $C_{i}=u_{-}-\dfrac{1}{2}\Big(\mp\sqrt{\epsilon^{2}_{2}+4\epsilon_{1}v_{-}^{2}}+
\epsilon_{2}\ln(\sqrt{\epsilon^{2}_{2}+4\epsilon_{1}v_{-}^{2}}\pm\epsilon_{2})\Big)$, $i=1,2$.
\\

In order to depict the geometric properties of rarefaction wave curves, we show the following Lemma.

\vspace{0.3cm}
\noindent{\bf Lemma 4.1.} \ For the back and forward rarefaction waves based on the given left state
$(u_{-},v_{-})$, we have
$$R_{1}(u_{-},v_{-}): \ \ \dfrac{\mbox{d}u}{\mbox{d}v}<0, \ \ \dfrac{\mbox{d}^{2}u}{\mbox{d}v^{2}}<0, \ \  \lim\limits_{v\rightarrow 0^{+}}u=u_{\ast\ast},$$
$$ \ \ R_{2}(u_{-},v_{-}): \ \ \dfrac{\mbox{d}u}{\mbox{d}v}>0, \ \ \dfrac{\mbox{d}^{2}u}{\mbox{d}v^{2}}<0, \ \  \lim\limits_{v\rightarrow +\infty}u=+\infty,$$
where $u_{\ast\ast}=\dfrac{1}{2}(-\epsilon_{2}+\epsilon_{2}\ln2\epsilon_{2})+C_{1}$.

\vspace{0.3cm}
\noindent{\bf Proof.}  From the second equation of \eqref{eq4.4}, we can obtain that
\begin{align}
\dfrac{\mbox{d}u}{\mbox{d}v}=\dfrac{\epsilon_{2}-\sqrt{\epsilon^{2}_{2}+4\epsilon_{1}v^{2}}}{2v}<0
\end{align}
and
\begin{align}
\dfrac{\mbox{d}^{2}u}{\mbox{d}v^{2}}=\dfrac{-2\epsilon_{2}\sqrt{\epsilon^{2}_{2}+4\epsilon_{1}v^{2}}+2\epsilon^{2}_{2}}{4v^{2}\sqrt{\epsilon^{2}_{2}+4\epsilon_{1}v^{2}}}<0,
\end{align}
which implies that the backward rarefaction wave curve $R_{1}$ is monotonically decreasing
and concave. Similarly, the second equation of \eqref{eq4.5} yields that
\begin{align}
\frac{\mbox{d}u}{\mbox{d}v}=\frac{\epsilon_{2}+\sqrt{\epsilon^{2}_{2}+4\epsilon_{1}v^{2}}}{2v}>0
\end{align}
and
\begin{align}
\frac{\mbox{d}^{2}u}{\mbox{d}v^{2}}=\frac{-2\epsilon_{2}\sqrt{\epsilon^{2}_{2}+4\epsilon_{1}v^{2}}-2\epsilon^{2}_{2}}{4v^{2}\sqrt{\epsilon^{2}_{2}+4\epsilon_{1}v^{2}}}<0.
\end{align} Thus, the forward rarefaction wave curve $R_{2}$ is monotonically increasing and concave.

The asymptotic properties of $R_{1}$ and $R_{2}$ can be easily obtained from \eqref{eq4.8} and \eqref{eq4.9}. The proof is completed. \ \ \ $\square$
\\

For a bounded discontinuity at $\xi=\sigma^{\epsilon_{1}\epsilon_{2}}$,  the Rankine-Hugoniot condition for system \eqref{eq1.1} can be written as
\begin{align}\label{eq4.14}
\left\{\begin{array}{l}
 \sigma^{\epsilon_{1}\epsilon_{2}}[u]=\big[\frac{1}{2}u^{2}+\frac{1}{2}\epsilon_{1}v^{2}], \\[0.5cm]
\sigma^{\epsilon_{1}\epsilon_{2}}[v]=\big[uv-\epsilon_{2}v].
\end{array}\right.
\end{align}
 Eliminating $\sigma^{\epsilon_{1}\epsilon_{2}}$ in \eqref{eq4.14}, we have
\begin{align}
\frac{1}{2}(v_{-}+v_{+})\Big(\frac{u_{-}-u_{+}}{v_{-}-v_{+}}\Big)^{2}-\epsilon_{2}\frac{u_{-}-u_{+}}{v_{-}-v_{+}}-\frac{1}{2}\epsilon_{1}(v_{-}+v_{+})=0.
\end{align}
Then one can obtain that
\begin{align}\label{eq4.16}
\frac{u_{-}-u_{+}}{v_{-}-v_{+}}=\frac{\epsilon_{2}\pm\sqrt{\epsilon^{2}_{2}+4\epsilon_{1}(v_{-}+v_{+})^{2}}}{v_{-}+v_{+}}.
\end{align}

In view of the classical Lax entropy conditions, the propagation speed $\sigma^{\epsilon_{1}\epsilon_{2}}_{1}$
should satisfy
\begin{align}\label{eq4.17}
\sigma^{\epsilon_{1}\epsilon_{2}}_{1}<\lambda_{1}(u_{-},v_{-}), \ \ \ \ \  \lambda_{1}(u_{+},v_{+})<\sigma^{\epsilon_{1}\epsilon_{2}}_{1}<\lambda_{2}(u_{+},v_{+})
\end{align}
for the backward shock wave, and the propagation speed $\sigma^{\epsilon_{1}\epsilon_{2}}_{2}$
should satisfy
\begin{align}\label{eq4.18}
\lambda_{1}(u_{-},v_{-})<\sigma^{\epsilon_{1}\epsilon_{2}}_{2}<\lambda_{2}(u_{-},v_{-}), \ \ \ \ \  \lambda_{2}(u_{+},v_{+})<\sigma^{\epsilon_{1}\epsilon_{2}}_{2}
\end{align}
for the forward shock wave.

Furthermore, from the second equation of \eqref{eq4.14}, we have
\begin{align}  \label{eq4.19}
\sigma^{\epsilon_{1}\epsilon_{2}}=u_{-}+\frac{v_{+}(u_{-}-u_{+})}{v_{-}-v_{+}}-\epsilon_{2}=u_{+}+\frac{v_{-}(u_{-}-u_{+})}{v_{-}-v_{+}}-\epsilon_{2}.
\end{align}
Together with \eqref{eq4.17}-\eqref{eq4.19}, one has
\begin{align}\label{eq4.20}
\frac{-2\epsilon_{1}v^{3}_{+}}{\epsilon_{2}+\sqrt{\epsilon^{2}_{2}+4\epsilon_{1}v_{+}^{2}}}<\frac{v_{-}v_{+}(u_{-}-u_{+})}{v_{-}-v_{+}}
<\frac{-2\epsilon_{1}v^{3}_{-}}{\epsilon_{2}+\sqrt{\epsilon^{2}_{2}+4\epsilon_{1}v_{-}^{2}}}
\end{align}
and
\begin{align}\label{eq4.21}
\frac{-2\epsilon_{1}v^{3}_{+}}{\epsilon_{2}-\sqrt{\epsilon^{2}_{2}+4\epsilon_{1}v_{+}^{2}}}<\frac{v_{-}v_{+}(u_{-}-u_{+})}{v_{-}-v_{+}}
<\frac{-2\epsilon_{1}v^{3}_{-}}{\epsilon_{2}-\sqrt{\epsilon^{2}_{2}+4\epsilon_{1}v_{-}^{2}}}.
\end{align}
\eqref{eq4.20} means that $v_{+}>v_{-}$, $u_{-}>u_{+}$, and the minus sign is taken in \eqref{eq4.16} for backward shock wave.
 In contrast, \eqref{eq4.21} means that $v_{+}<v_{-}$, $u_{-}>u_{+}$, and the plus sign is chosen in \eqref{eq4.16} for forward shock wave.

Let us fix $(u_{-},v_{-})$ in the $(u,v)$ phase plane,
the state $(u, v)$ which is connected to $(u_{-},v_{-})$ on the right by shock waves are given below.
The backward shock wave curve in the phase plane is
\begin{align}\label{eq4.22}
S_{1}(u_{-},v_{-}):\left\{\begin{array}{l}
 \sigma^{\epsilon_{1}\epsilon_{2}}_{1} =u_{-}+\dfrac{v(u-u_{-})}{v-v_{-}}-\epsilon_{2}, \cr\noalign {\vskip2truemm}
\dfrac{u-u_{-}}{v-v_{-}}=\dfrac{\epsilon_{2}-\sqrt{\epsilon^{2}_{2}+4\epsilon_{1}(v+v_{-})^{2}}}{v+v_{-}},\cr\noalign {\vskip2truemm}
u<u_{-},\ \ v>v_{-}
\end{array}\right.
\end{align}
and the forward shock wave curve in the phase plane is
\begin{align}\label{eq4.23}
S_{2}(u_{-},v_{-}):\left\{\begin{array}{l}
 \sigma^{\epsilon_{1}\epsilon_{2}}_{2} =u+\dfrac{v_{-}(u-u_{-})}{v-v_{-}}-\epsilon_{2}, \cr\noalign {\vskip2truemm}
\dfrac{u-u_{-}}{v-v_{-}}=\dfrac{\epsilon_{2}+\sqrt{\epsilon^{2}_{2}+4\epsilon_{1}(v+v_{-})^{2}}}{v+v_{-}},\cr\noalign {\vskip2truemm}
 u<u_{-},\ \ v<v_{-}.
\end{array}\right.
\end{align}

The shock wave curves possess the following geometric properties.

\vspace{0.3cm}
\noindent{\bf Lemma 4.2.} \ For the back and forward shock waves based on the given left state
$(u_{-},v_{-})$, we have
$$S_{1}(u_{-},v_{-}): \ \ \dfrac{\mbox{d}u}{\mbox{d}v}<0, \ \ \dfrac{\mbox{d}^{2}u}{\mbox{d}v^{2}}<0, \ \  \lim\limits_{v\rightarrow +\infty}u=-\infty,$$
$$ \ \ S_{2}(u_{-},v_{-}): \ \ \dfrac{\mbox{d}u}{\mbox{d}v}>0, \ \ \dfrac{\mbox{d}^{2}u}{\mbox{d}v^{2}}<0, \ \  \lim\limits_{v\rightarrow 0^{+}}u=u_{\ast},$$
where $u_{\ast}=u_{-}-(\epsilon_{2}+\sqrt{\epsilon^{2}_{2}+4\epsilon_{1}v^{2}_{-}})$.

\vspace{0.3cm}
\noindent{\bf Proof.} By a simple calculation, it follows from the second equation of \eqref{eq4.22} that
\begin{align}
\dfrac{\mbox{d}u}{\mbox{d}v}=\dfrac{\epsilon_{2}-\sqrt{\epsilon^{2}_{2}+4\epsilon_{1}(v+v_{-})^{2}}}{v+v_{-}}
+(v-v_{-})\dfrac{\epsilon^{2}_{2}-\epsilon_{2}\sqrt{\epsilon^{2}_{2}+4\epsilon_{1}(v+v_{-})^{2}}}{(v+v_{-})^{2}\sqrt{\epsilon^{2}_{2}+4\epsilon_{1}(v+v_{-})^{2}}}<0,
\end{align}
\begin{align}
\begin{array}{l}
\dfrac{\mbox{d}^{2}u}{\mbox{d}v^{2}}=\dfrac{4v_{-}(v+v_{-})\big(\epsilon^{2}_{2}+4\epsilon_{1}(v+v_{-})^{2}\big)\Big(\epsilon^{2}_{2}-\epsilon_{2}\sqrt{\epsilon^{2}_{2}+4\epsilon_{1}(v+v_{-})^{2}}\Big)}
{(v+v_{-})^{4}(\epsilon^{2}_{2}+4\epsilon_{1}(v+v_{-})^{2})^{\frac{3}{2}}}\cr\noalign {\vskip3truemm}
\ \ \ \ \ \ \ \ \ +\dfrac{-4\epsilon_{1}\epsilon_{2}(v+v_{-})^{3}(v-v_{-})}
{(v+v_{-})^{4}(\epsilon^{2}_{2}+4\epsilon_{1}(v+v_{-})^{2})^{\frac{3}{2}}}<0.
\end{array}
\end{align}
Therefore, the backward shock wave curve $S_{1}$ is monotonically decreasing and concave.

Similarly, from the second equation of \eqref{eq4.23}, we can get
\begin{align}
\dfrac{\mbox{d}u}{\mbox{d}v}=\dfrac{\epsilon_{2}+\sqrt{\epsilon^{2}_{2}+4\epsilon_{1}(v+v_{-})^{2}}}{v+v_{-}}
+(v-v_{-})\dfrac{-\epsilon^{2}_{2}-\epsilon_{2}\sqrt{\epsilon^{2}_{2}+4\epsilon_{1}(v+v_{-})^{2}}}{(v+v_{-})^{2}\sqrt{\epsilon^{2}_{2}+4\epsilon_{1}(v+v_{-})^{2}}}>0,
\end{align}
\begin{align}
\begin{array}{l}
\dfrac{\mbox{d}^{2}u}{\mbox{d}v^{2}}=\dfrac{4v_{-}(v+v_{-})\big(\epsilon^{2}_{2}+4\epsilon_{1}(v+v_{-})^{2}\big)\Big(-\epsilon^{2}_{2}-\epsilon_{2}\sqrt{\epsilon^{2}_{2}+4\epsilon_{1}(v+v_{-})^{2}}\Big)}
{(v+v_{-})^{4}(\epsilon^{2}_{2}+4\epsilon_{1}(v+v_{-})^{2})^{\frac{3}{2}}}\cr\noalign {\vskip3truemm}
\ \ \ \ \ \ \ \ \ +\dfrac{4\epsilon_{1}\epsilon_{2}(v+v_{-})^{3}(v-v_{-})}
{(v+v_{-})^{4}(\epsilon^{2}_{2}+4\epsilon_{1}(v+v_{-})^{2})^{\frac{3}{2}}}<0.
\end{array}
\end{align}
That is, the forward shock curve $S_{2}$ is monotonically increasing and concave.
In addition, the asymptotic properties of $S_{1}$ and $S_{2}$ are obviously checked from  \eqref{eq4.22} and \eqref{eq4.23}. We finish the proof. \ \ \ $\square$\\

Through the above analysis, for a given left state $(u_{-},v_{-})$,
we know that the elementary wave curves divide the half-upper $(u,v)$ phase plane into four regions.
The structures of  Riemann solutions for \eqref{eq1.5} and \eqref{eq1.6} depend on the position of the right state $(u_{+},v_{+})$.
The solution is $R_{1}+R_{2}$ when $(u_{+},v_{+})\in R_{1}R_{2}(u_{-},v_{-})$,
$S_{1}+R_{2}$ when $(u_{+},v_{+})\in S_{1}R_{2}(u_{-},v_{-}) $, $R_{1}+S_{2}$ when $(u_{+},v_{+})\in  R_{1}S_{2}(u_{-},v_{-}) $, and $S_{1}+S_{2}$ when $(u_{+},v_{+})\in S_{1}S_{2}(u_{-},v_{-})$(see
Fig.2.).

\section{Formation of delta shock wave and vacuum when $\epsilon_{1},\epsilon_{2}\rightarrow0$}

In this section, we investigate the formation process of delta shock wave and vacuum
in the Riemann solutions to \eqref{eq1.1} and \eqref{eq1.6} as $\epsilon_{1},\epsilon_{2}\rightarrow0$.
Our discussions will be divided into the following two cases.

\subsection{Formation of delta shock wave for the system \eqref{eq1.1}}

We at first discuss the limit behavior of the Riemann solutions as $\epsilon_{1},\epsilon_{2}\rightarrow0$.

When $(u_{+},v_{+})\in S_{1}S_{2}(u_{-},v_{-})$, the solution of Riemann problem \eqref{eq1.1} and \eqref{eq1.6} consists of
$S_{1}$ and $S_{2}$. Denoted the intermediate state
by $(u^{\epsilon_{1}\epsilon_{2}}_{\ast},v^{\epsilon_{1}\epsilon_{2}}_{\ast})$, then we have
\begin{align}\label{eq5.1}
S_{1}(u_{-},v_{-}):\left\{\begin{array}{l}
 \sigma_{1}^{\epsilon_{1}\epsilon_{2}} =u_{-}+\frac{v^{\epsilon_{1}\epsilon_{2}}_{\ast}(u^{\epsilon_{1}\epsilon_{2}}_{\ast}-u_{-})}{v^{\epsilon_{1}\epsilon_{2}}_{\ast}-v_{-}}-\epsilon_{2}, \cr\noalign {\vskip2truemm}
u^{\epsilon_{1}\epsilon_{2}}_{\ast}=u_{-}+(v^{\epsilon_{1}\epsilon_{2}}_{\ast}-v_{-})\dfrac{\epsilon_{2}-\sqrt{\epsilon^{2}_{2}+4\epsilon_{1}(v^{\epsilon_{1}\epsilon_{2}}_{\ast}+v_{-})^{2}}}{v^{\epsilon_{1}\epsilon_{2}}_{\ast}+v_{-}}, \cr\noalign {\vskip2truemm}
 u^{\epsilon_{1}\epsilon_{2}}_{\ast}<u_{-}, \ \ v^{\epsilon_{1}\epsilon_{2}}_{\ast}>v_{-}
\end{array}\right.
\end{align}
and
\begin{align}\label{eq5.2}
S_{2}(u_{-},v_{-}):\left\{\begin{array}{l}
 \sigma^{\epsilon_{1}\epsilon_{2}}_{2} =u_{+}+\frac{v^{\epsilon_{1}\epsilon_{2}}_{\ast}(u_{+}-u^{\epsilon_{1}\epsilon_{2}}_{\ast})}{v_{+}-v^{\epsilon_{1}\epsilon_{2}}_{\ast}}-\epsilon_{2}, \cr\noalign {\vskip2truemm}
u_{+}=u^{\epsilon_{1}\epsilon_{2}}_{\ast}+(v_{+}-v^{\epsilon_{1}\epsilon_{2}}_{\ast})\dfrac{\epsilon_{2}+\sqrt{\epsilon^{2}_{2}+4\epsilon_{1}(v^{\epsilon_{1}\epsilon_{2}}_{\ast}+v_{+})^{2}}}{v^{\epsilon_{1}\epsilon_{2}}_{\ast}+v_{+}}, \cr\noalign {\vskip2truemm}
 u^{\epsilon_{1}\epsilon_{2}}_{\ast}>u_{+}, \ \ v^{\epsilon_{1}\epsilon_{2}}_{\ast}>v_{+}.
\end{array}\right.
\end{align}

We need to establish the following lemmas.

\vspace{0.3cm}
\noindent{ \bf Lemma 5.1} $\lim\limits_{\epsilon_{1},\epsilon_{2}\rightarrow0}v^{\epsilon_{1}\epsilon_{2}}_{\ast}=+\infty$.

\vspace{0.3cm}
\noindent{ \bf Proof.} From the second equations of \eqref{eq5.1} and \eqref{eq5.2} yields that
\begin{align}\label{eq5.3}
\begin{array}{l}
u_{-}-u_{+}=(v^{\epsilon_{1}\epsilon_{2}}_{\ast}-v_{+})\dfrac{\epsilon_{2}+\sqrt{\epsilon^{2}_{2}+4\epsilon_{1}(v^{\epsilon_{1}\epsilon_{2}}_{\ast}+v_{+})^{2}}}{v^{\epsilon_{1}\epsilon_{2}}_{\ast}+v_{+}}\cr\noalign {\vskip2truemm}
\hspace{2cm}+(v_{-}-v^{\epsilon_{1}\epsilon_{2}}_{\ast})\dfrac{\epsilon_{2}-\sqrt{\epsilon^{2}_{2}+4\epsilon_{1}(v^{\epsilon_{1}\epsilon_{2}}_{\ast}+v_{-})^{2}}}{v^{\epsilon_{1}\epsilon_{2}}_{\ast}+v_{-}}.
\end{array}
\end{align}
If $\lim\limits_{\epsilon_{1},\epsilon_{2}\rightarrow0}v^{\epsilon_{1}\epsilon_{2}}_{\ast}=M\in(\mbox{max}({v_{-},v_{+}}),+\infty)$,
then by passing to the limit of \eqref{eq5.3} as $\epsilon_{1},\epsilon_{2}\rightarrow0$,
we directly get $u_{-}-u_{+}=0$, which is a contradiction. Thus, one has
$\lim\limits_{\epsilon_{1},\epsilon_{2}\rightarrow0}v^{\epsilon_{1}\epsilon_{2}}_{\ast}=+\infty$. This concludes the proof.  \ \ \ $\square$
\\

Lemma 5.1 shows that  the intermediate state $v^{\epsilon_{1}\epsilon_{2}}_{\ast}$ between the two shock waves
becomes singular
when $\epsilon_{1}$ and $\epsilon_{2}$ drop to 0.

\vspace{0.3cm}
\noindent{ \bf Lemma 5.2} $\lim\limits_{\epsilon_{1},\epsilon_{2}\rightarrow0}2\sqrt{\epsilon_{1}}v^{\epsilon_{1}\epsilon_{2}}_{\ast}=\dfrac{u_{-}-u_{+}}{2}$.

\vspace{0.3cm}
\noindent{ \bf Proof.} From Lemma 5.1 and \eqref{eq5.3}, we have
\begin{align}
\begin{array}{l}
u_{-}-u_{+}=\lim\limits_{\epsilon_{1},\epsilon_{2}\rightarrow0}\big(\sqrt{4\epsilon_{1}(v^{\epsilon_{1}\epsilon_{2}}_{\ast}+v_{+})^{2}}
+\sqrt{4\epsilon_{1}(v^{\epsilon_{1}\epsilon_{2}}_{\ast}+v_{-})^{2}}\big)\cr\noalign {\vskip2truemm}
\hspace{1.5cm}=\lim\limits_{\epsilon_{1},\epsilon_{2}\rightarrow0}2\sqrt{\epsilon_{1}}(2v^{\epsilon_{1}\epsilon_{2}}_{\ast}+v_{+}+v_{-}).
\end{array}
\end{align}
Then, the desired result can be immediately obtained. The proof is completed.  \ \ \ $\square$

\vspace{0.3cm}
\noindent{ \bf Lemma 5.3} $ \lim\limits_{\epsilon_{1},\epsilon_{2}\rightarrow0}u^{\epsilon_{1}\epsilon_{2}}_{\ast}= \lim\limits_{\epsilon_{1},\epsilon_{2}\rightarrow0}\sigma^{\epsilon_{1}\epsilon_{2}}_{1}=\lim\limits_{\epsilon_{1},\epsilon_{2}\rightarrow0}\sigma^{\epsilon_{1}\epsilon_{2}}_{2}
=\dfrac{u_{-}+u_{+}}{2}.$

\vspace{0.2cm}
\noindent{ \bf Proof.}
Combining \eqref{eq5.1} and \eqref{eq5.2} with Lemma 5.2, we obtain that
\begin{align}
\lim\limits_{\epsilon_{1},\epsilon_{2}\rightarrow0}u^{\epsilon_{1}\epsilon_{2}}_{\ast}=\lim\limits_{\epsilon_{1},\epsilon_{2}\rightarrow0}\Big(u_{-}-\sqrt{\epsilon^{2}_{2}+4\epsilon_{1}(v^{\epsilon_{1}\epsilon_{2}}_{\ast}+v_{-})^{2}}\Big)
=\frac{u_{-}+u_{+}}{2}
\end{align}
and
\begin{align}
\lim\limits_{\epsilon_{1},\epsilon_{2}\rightarrow0}\sigma^{\epsilon_{1}\epsilon_{2}}_{1}=\lim\limits_{\epsilon_{1},\epsilon_{2}\rightarrow0}\Big(u_{-}+v^{\epsilon_{1}\epsilon_{2}}_{\ast}\dfrac{\epsilon_{2}-\sqrt{\epsilon^{2}_{2}+4\epsilon_{1}(v^{\epsilon_{1}\epsilon_{2}}_{\ast}+v_{-})^{2}}}{v^{\epsilon_{1}\epsilon_{2}}_{\ast}+v_{-}}-\epsilon_{2}\Big)
=\frac{u_{-}+u_{+}}{2}.
\end{align}
Similarly, we can prove that $\lim\limits_{\epsilon_{2},\epsilon_{2}\rightarrow0}\sigma^{\epsilon_{1}\epsilon_{2}}_{2}=\dfrac{u_{-}+u_{+}}{2}$.
Thus, this lemma is right.  \ \ \ $\square$
\\

Lemma 5.3 means that, when $\epsilon_{1},\epsilon_{2}\rightarrow0$, $S_{1}$ and $S_{2}$ coincide, and
the velocities $\sigma^{\epsilon_{1}\epsilon_{2}}_{1}$ and $\sigma^{\epsilon_{1}\epsilon_{2}}_{2}$
approach to the quantity $\sigma$ given in \eqref{eq2.8}, which is just the propagating speed of the delta shock wave of the
transport equations.

\vspace{0.3cm}
\noindent{ \bf Lemma 5.4} $ \lim\limits_{\epsilon_{1},\epsilon_{2}\rightarrow0}\dis\int_{\sigma^{\epsilon_{1}\epsilon_{2}}_{1}}^{\sigma^{\epsilon_{1}\epsilon_{2}}_{2}}v^{\epsilon_{1}\epsilon_{2}}_{\ast}\mbox{d}\xi
=\sigma[v]-[uv]=\frac{1}{2}(v_{-}+v_{+})(u_{-}-u_{+}).$

\vspace{0.3cm}
\noindent{ \bf Proof.} $S_{1}$ and $S_{2}$ satisfy the following Rankine-Hugoniot relation
\begin{align}
\left\{\begin{array}{l}
\sigma^{\epsilon_{1}\epsilon_{2}}_{1}(v_{-}-v^{\epsilon_{1}\epsilon_{2}}_{\ast})=u_{-}v_{-}-\epsilon_{2}v_{-}-u^{\epsilon_{1}\epsilon_{2}}_{\ast}v^{\epsilon_{1}\epsilon_{2}}_{\ast}+\epsilon_{2}v^{\epsilon_{1}\epsilon_{2}}_{\ast}, \cr\noalign {\vskip2truemm}
\sigma^{\epsilon_{1}\epsilon_{2}}_{2}(v^{\epsilon_{1}\epsilon_{2}}_{\ast}-v_{+})=u^{\epsilon_{1}\epsilon_{2}}_{\ast}v^{\epsilon_{1}\epsilon_{2}}_{\ast}-\epsilon_{2}v^{\epsilon_{1}\epsilon_{2}}_{\ast}-u_{+}v_{+}+\epsilon_{2}v_{+}.
\end{array}\right.
\end{align}
Adding them together, one has
$$
\begin{array}{l}
\lim\limits_{\epsilon_{1},\epsilon_{2}\rightarrow0}v^{\epsilon_{1}\epsilon_{2}}_{\ast}(\sigma^{\epsilon_{1}\epsilon_{2}}_{2}-\sigma^{\epsilon_{1}\epsilon_{2}}_{1})=\lim\limits_{\epsilon_{1},\epsilon_{2}\rightarrow0}(\sigma^{\epsilon_{1}\epsilon_{2}}_{2}v_{+}-\sigma^{\epsilon_{1}\epsilon_{2}}_{1}v_{-}+u_{-}v_{-}-u_{+}v_{+}+\epsilon_{2}v_{+}-\epsilon_{2}v_{-})\cr\noalign {\vskip3truemm}
\hspace{4.3cm}=\sigma[v]-[uv]\cr\noalign {\vskip2truemm}
\hspace{4.3cm}=\dfrac{1}{2}(v_{-}+v_{+})(u_{-}-u_{+}).
\end{array}
$$
This result is consistent with the strength of delta shock wave given by \eqref{eq2.8} in Section 2.
Thus, we complete the proof of this lemma.  \ \ \ $\square$
\\

Now, we present the limit of solutions to  \eqref{eq1.1} and  \eqref{eq1.6}
as two-parameter flux perturbation vanishes.
The following theorem is necessary because it gives a very nice depiction of the limit.

\vspace{0.3cm}
\noindent{ \bf Theorem 5.1.} Let $u_{-}>u_{+}$. For any fixed $\epsilon_{1}$, $\epsilon_{2}>0$, assume that $(u^{\epsilon_{1}\epsilon_{2}}, v^{\epsilon_{1}\epsilon_{2}})$ is a two-shock
Riemann solution of the system \eqref{eq1.1} and \eqref{eq1.6} constructed in Section 4. Then, as
$\epsilon_{1},\epsilon_{2}\rightarrow0$, the limit of solution $(u^{\epsilon_{1}\epsilon_{2}}, v^{\epsilon_{1}\epsilon_{2}})$
is a delta shock wave of  \eqref{eq1.2} and \eqref{eq1.6}  connecting two constant states $(u_{\pm},v_{\pm})$, which is expressed by \eqref{eq2.6} with \eqref{eq2.8}.

\vspace{0.3cm}
\noindent{\bf Proof.}\ \ (1). Firstly, set $\xi=x/t$. For any $\epsilon_{1}$, $\epsilon_{2}>0$, the two-shock Riemann solution
can be expressed as
\begin{align}\label{eq5.8}
(u^{\epsilon_{1}\epsilon_{2}},v^{\epsilon_{1}\epsilon_{2}})(\xi)=\left\{\begin{array}{l}
 (u_{-},v_{-}),  \ \ \ \ \ \ \ \ \ \xi<\sigma^{\epsilon_{1}\epsilon_{2}}_{1},\cr\noalign {\vskip2truemm}
(u^{\epsilon_{1}\epsilon_{2}}_{\ast},v^{\epsilon_{1}\epsilon_{2}}_{\ast}), \ \ \ \ \sigma^{\epsilon_{1}\epsilon_{2}}_{1}<\xi<\sigma^{\epsilon_{1}\epsilon_{2}}_{2},\cr\noalign {\vskip2truemm}
 (u_{+},v_{+}),  \ \ \ \ \ \ \ \ \  \xi>\sigma^{\epsilon_{1}\epsilon_{2}}_{2}
\end{array}\right.
\end{align}
which satisfies the weak formulations
\begin{align}\label{eq5.9}
\int_{-\infty}^{+\infty}(\xi u^{\epsilon_{1}\epsilon_{2}}-\frac{1}{2}(u^{\epsilon_{1}\epsilon_{2}})^{2}-\frac{1}{2}\epsilon_{1}(v^{\epsilon_{1}\epsilon_{2}})^{2})\psi'd\xi+\int_{-\infty}^{+\infty}u^{\epsilon_{1}\epsilon_{2}}\psi d\xi=0,
\end{align}
and
\begin{align}\label{eq5.10}
\int_{-\infty}^{+\infty}(\xi v^{\epsilon_{1}\epsilon_{2}}-(u^{\epsilon_{1}\epsilon_{2}}v^{\epsilon_{1}\epsilon_{2}}-\epsilon_{2}v^{\epsilon_{1}\epsilon_{2}}))\psi'd\xi+\int_{-\infty}^{+\infty}v^{\epsilon_{1}\epsilon_{2}}\psi d\xi=0
\end{align}
for any function $\psi\in C_{0}^{1}(-\infty, +\infty)$.

(2). Secondly, we consider the limits of $u^{\epsilon_{1}\epsilon_{2}}$ and $v^{\epsilon_{1}\epsilon_{2}}$ depending on $\xi$.
Dividing the integral interval $(-\infty,+\infty)$ into $(-\infty,\sigma^{\epsilon_{1}\epsilon_{2}}_{1})$,
$(\sigma^{\epsilon_{1}\epsilon_{2}}_{1},\sigma^{\epsilon_{1}\epsilon_{2}}_{2})$ and $(\sigma^{\epsilon_{1}\epsilon_{2}}_{2},+\infty)$,
for the first integral on the left side of \eqref{eq5.10}, we have
\begin{align}
\begin{array}{l}
\dis\int_{-\infty}^{+\infty}(u^{\epsilon_{1}\epsilon_{2}}v^{\epsilon_{1}\epsilon_{2}}-\epsilon_{2}v^{\epsilon_{1}\epsilon_{2}}-\xi v^{\epsilon_{1}\epsilon_{2}})\psi'd\xi\cr\noalign {\vskip3truemm}
=\Big(\dis\int_{-\infty}^{\sigma^{\epsilon_{1}\epsilon_{2}}_{1}}+\dis\int_{\sigma^{\epsilon_{1}\epsilon_{2}}_{1}}^{\sigma^{\epsilon_{1}\epsilon_{2}}_{2}}+\dis\int_{\sigma^{\epsilon_{1}\epsilon_{2}}_{2}}^{+\infty}\Big)(u^{\epsilon_{1}\epsilon_{2}}v^{\epsilon_{1}\epsilon_{2}}-\epsilon_{2}v^{\epsilon_{1}\epsilon_{2}}-\xi v^{\epsilon_{1}\epsilon_{2}})\psi'd\xi.
\end{array}
\end{align}

By integrating by parts and using Lemmas 5.3 and 5.4, one can calculate that
\begin{align}
\begin{array}{l}
\lim\limits_{\epsilon_{1},\epsilon_{2}\rightarrow0}\Big(\dis \int_{-\infty}^{\sigma^{\epsilon_{1}\epsilon_{2}}_{1}}+\dis\int_{\sigma^{\epsilon_{1}\epsilon_{2}}_{2}}^{+\infty}\Big)
(u^{\epsilon_{1}\epsilon_{2}}v^{\epsilon_{1}\epsilon_{2}}-\epsilon_{2}v^{\epsilon_{1}\epsilon_{2}}-\xi v^{\epsilon_{1}\epsilon_{2}})\psi'd\xi\cr\noalign {\vskip3truemm}
= \lim\limits_{\epsilon_{1},\epsilon_{2}\rightarrow0}\dis\int_{-\infty}^{\sigma^{\epsilon_{1}\epsilon_{2}}_{1}}(u_{-}v_{-}-\epsilon_{2}v_{-}-\xi v_{-})\psi'd\xi
+\lim\limits_{\epsilon_{1},\epsilon_{2}\rightarrow0}\dis\int_{\sigma^{\epsilon_{1}\epsilon_{2}}_{2}}^{+\infty}(u_{+}v_{+}-\epsilon_{2}v_{+}-\xi v_{+})\psi'd\xi\cr\noalign {\vskip3truemm}
=(\sigma[v]-[uv])\psi(\sigma)+\dis\int_{-\infty}^{+\infty}\psi H_{v}(\xi-\sigma)d\xi,
\end{array}
\end{align}
where
$$
H_{v}(x)=\left\{\begin{array}{ll}
v_-,&x<0,\\
v_+,&x>0.
\end{array}\right.
$$
Meanwhile, we have
\begin{align}\label{eq5.13}
\begin{array}{l}
\lim\limits_{\epsilon_{1},\epsilon_{2}\rightarrow0}\dis \int_{\sigma^{\epsilon_{1}\epsilon_{2}}_{1}}^{\sigma^{\epsilon_{1}\epsilon_{2}}_{2}}(u^{\epsilon_{1}\epsilon_{2}}v^{\epsilon_{1}\epsilon_{2}}-\epsilon_{2}v^{\epsilon_{1}\epsilon_{2}}-\xi v^{\epsilon_{1}\epsilon_{2}})\psi'd\xi \cr\noalign{\vskip3truemm}
=\lim\limits_{\epsilon_{1},\epsilon_{2}\rightarrow0}\Big(u^{\epsilon_{1}\epsilon_{2}}_{\ast}v^{\epsilon_{1}\epsilon_{2}}_{\ast}(\psi(\sigma^{\epsilon_{1}\epsilon_{2}}_{2})-\psi(\sigma^{\epsilon_{1}\epsilon_{2}}_{1}))-\epsilon_{2}v^{\epsilon_{1}\epsilon_{2}}_{\ast}\psi(\sigma^{\epsilon_{1}\epsilon_{2}}_{2})+\epsilon_{2}v^{\epsilon_{1}\epsilon_{2}}_{\ast}\psi(\sigma^{\epsilon_{1}\epsilon_{2}}_{1})
\cr\noalign {\vskip3truemm}
\ \ \ \ -\sigma^{\epsilon_{1}\epsilon_{2}}_{2}v^{\epsilon_{1}\epsilon_{2}}_{\ast}\psi(\sigma^{\epsilon_{1}\epsilon_{2}}_{2})+\sigma^{\epsilon_{1}\epsilon_{2}}_{1}v^{\epsilon_{1}\epsilon_{2}}_{\ast}\psi(\sigma^{\epsilon_{1}\epsilon_{2}}_{1})\Big)
+\lim\limits_{\epsilon_{1},\epsilon_{2}\rightarrow0}\dis\int_{\sigma^{\epsilon_{1}\epsilon_{2}}_{1}}^{\sigma^{\epsilon_{1}\epsilon_{2}}_{2}}v^{\epsilon_{1}\epsilon_{2}}_{\ast}\psi d\xi.
\end{array}
\end{align}
Considering Lemmas 5.3 and 5.4, we can obtain
\begin{align}
\begin{array}{l}
\lim\limits_{\epsilon_{1},\epsilon_{2}\rightarrow0}u^{\epsilon_{1}\epsilon_{2}}_{\ast}v^{\epsilon_{1}\epsilon_{2}}_{\ast}(\psi(\sigma^{\epsilon_{1}\epsilon_{2}}_{2})-\psi(\sigma^{\epsilon_{1}\epsilon_{2}}_{1}))=\lim\limits_{\epsilon_{1},\epsilon_{2}\rightarrow0}
u^{\epsilon_{1}\epsilon_{2}}_{\ast}v^{\epsilon_{1}\epsilon_{2}}_{\ast}(\sigma^{\epsilon_{1}\epsilon_{2}}_{2}-\sigma^{\epsilon_{1}\epsilon_{2}}_{1})
\dfrac{\psi(\sigma^{\epsilon_{1}\epsilon_{2}}_{2})-\psi(\sigma^{\epsilon_{1}\epsilon_{2}}_{1})}{\sigma^{\epsilon_{1}\epsilon_{2}}_{2}-\sigma^{\epsilon_{1}\epsilon_{2}}_{1}}\cr\noalign {\vskip3truemm}
\hspace{6.2cm}=(\sigma[v]-[uv])\sigma\psi'(\sigma).
\end{array}
\end{align}

Similarly, we can obtain that
\begin{align}
\begin{array}{l}
\lim\limits_{\epsilon_{1},\epsilon_{2}\rightarrow0}\epsilon_{2}v^{\epsilon_{1}\epsilon_{2}}_{\ast}\psi(\sigma^{\epsilon_{1}\epsilon_{2}}_{1})-\epsilon_{2}v^{\epsilon_{1}\epsilon_{2}}_{\ast}\psi(\sigma^{\epsilon_{1}\epsilon_{2}}_{2})\cr\noalign {\vskip3truemm}
 =\lim\limits_{\epsilon_{1},\epsilon_{2}\rightarrow0}\epsilon_{2}v^{\epsilon_{1}\epsilon_{2}}_{\ast}(\sigma^{\epsilon_{1}\epsilon_{2}}_{1}-\sigma^{\epsilon_{1}\epsilon_{2}}_{2})\dfrac{\psi(\sigma^{\epsilon_{1}\epsilon_{2}}_{1})-\psi(\sigma^{\epsilon_{1}\epsilon_{2}}_{2})}
{\sigma^{\epsilon_{1}\epsilon_{2}}_{1}-\sigma^{\epsilon_{1}\epsilon_{2}}_{2}}\cr\noalign {\vskip3truemm}
=0
\end{array}
\end{align}
and
\begin{align}
\begin{array}{l}
\lim\limits_{\epsilon_{1},\epsilon_{2}\rightarrow0}\sigma^{\epsilon_{1}\epsilon_{2}}_{1}v^{\epsilon_{1}\epsilon_{2}}_{\ast}\psi(\sigma^{\epsilon_{1}\epsilon_{2}}_{1})-\sigma^{\epsilon_{1}\epsilon_{2}}_{2}v^{\epsilon_{1}\epsilon_{2}}_{\ast}\psi(\sigma^{\epsilon_{1}\epsilon_{2}}_{2})\cr\noalign {\vskip3truemm}
 =\lim\limits_{\epsilon_{1},\epsilon_{2}\rightarrow0}v^{\epsilon_{1}\epsilon_{2}}_{\ast}(\sigma^{\epsilon_{1}\epsilon_{2}}_{1}-\sigma^{\epsilon_{1}\epsilon_{2}}_{2})\dfrac{\sigma_{1}^{\epsilon_{1}\epsilon_{2}}\psi(\sigma^{\epsilon_{1}\epsilon_{2}}_{1})-\sigma^{\epsilon_{1}\epsilon_{2}}_{2}\psi(\sigma^{\epsilon_{1}\epsilon_{2}}_{2})}{\sigma^{\epsilon_{1}\epsilon_{2}}_{1}-\sigma^{\epsilon_{1}\epsilon_{2}}_{2}}\cr\noalign
{\vskip3truemm}
=-(\sigma[v]-[uv])(\sigma\psi'(\sigma)+\psi(\sigma)).
\end{array}
\end{align}

Moreover, the Lemma 5.4 yields that
\begin{align}
\begin{array}{l}
\lim\limits_{\epsilon_{1},\epsilon_{2}\rightarrow0}\dis\int_{\sigma^{\epsilon_{1}\epsilon_{2}}_{1}}^{\sigma^{\epsilon_{1}\epsilon_{2}}_{2}}v^{\epsilon_{1}\epsilon_{2}}_{\ast}\psi\mbox{d}\xi
=\lim\limits_{\epsilon_{1},\epsilon_{2}\rightarrow0}(\sigma^{\epsilon_{1}\epsilon_{2}}_{2}-\sigma^{\epsilon_{1}\epsilon_{2}}_{1})v^{\epsilon_{1}\epsilon_{2}}_{\ast}
\cdot \lim\limits_{\epsilon_{1},\epsilon_{2}\rightarrow0}\frac{1}{\sigma^{\epsilon_{1}\epsilon_{2}}_{2}-\sigma^{\epsilon_{1}\epsilon_{2}}_{1}}
\dis\int_{\sigma^{\epsilon_{1}\epsilon_{2}}_{1}}^{\sigma^{\epsilon_{1}\epsilon_{2}}_{2}}\psi(\xi)\mbox{d}\xi\cr\noalign {\vskip3truemm}
\hspace{3.8cm}=(\sigma[v]-[uv])\psi(\sigma).
\end{array}
\end{align}
Returning to \eqref{eq5.13}, we immediately get
\begin{align}
\lim\limits_{\epsilon_{1},\epsilon_{2}\rightarrow0}\dis \int_{\sigma^{\epsilon_{1}\epsilon_{2}}_{1}}^{\sigma^{\epsilon_{1}\epsilon_{2}}_{2}}(u^{\epsilon_{1}\epsilon_{2}}v^{\epsilon_{1}\epsilon_{2}}-\epsilon_{2}v^{\epsilon_{1}\epsilon_{2}}-\xi v^{\epsilon_{1}\epsilon_{2}})\psi'd\xi=0.
\end{align}
Therefore
\begin{align}
\lim\limits_{\epsilon_{1},\epsilon_{2}\rightarrow0}\dis\int_{-\infty}^{+\infty}(v^{\epsilon_{1}\epsilon_{2}}-H_{v}(\xi-\sigma))\psi d\xi=(\sigma[v]-[uv])\psi(\sigma).
\end{align}

Now, let us focus on \eqref{eq5.9}. Noticing the fact that  $\epsilon_{1}(v^{\epsilon_{1}\epsilon_{2}}_{\ast})^{2}$ is bounded
and using the same way as above, one can obtain that
\begin{align}\label{eq5.20}
\begin{array}{l}
\hspace{-3cm}\lim\limits_{\epsilon_{1},\epsilon_{2}\rightarrow0}\dis\int_{-\infty}^{+\infty}(u^{\epsilon_{1}\epsilon_{2}}-H_{u}(\xi-\sigma))\psi d\xi=(\sigma[u]-\big[\frac{1}{2}u^{2}\big])\psi(\sigma)\cr\noalign {\vskip3truemm}
\hspace{3.1cm}=\Big(\dfrac{u_{-}+u_{+}}{2}[u]-\big[\frac{1}{2}u^{2}\big]\Big)\psi(\sigma),\cr\noalign {\vskip3truemm}
\hspace{3.1cm}=0,
\end{array}
\end{align}
where
$$
H_{u}(x)=\left\{\begin{array}{ll}
u_-,&x<0,\\
u_+,&x>0.
\end{array}\right.
$$

(3). Finally, we study the limits of $u^{\epsilon_{1}\epsilon_{2}}(x,t)$ and $v^{\epsilon_{1}\epsilon_{2}}(x,t)$ depending on $t$. For any test function $\phi(x,t)\in C_{0}^{+\infty}(R\times R^{+})$,
 by \eqref{eq5.20}, we have
\begin{align}
\begin{array}{l}
\lim\limits_{\epsilon_{1},\epsilon_{2}\rightarrow0}\dis \int_{0}^{+\infty} \int_{-\infty}^{+\infty}u^{\epsilon_{1}\epsilon_{2}}(x/t)\phi(x,t)
dxdt\cr\noalign{\vskip3truemm}
=\lim\limits_{\epsilon_{1},\epsilon_{2}\rightarrow0}\dis \int_{0}^{+\infty} t(\int_{-\infty}^{+\infty}u^{\epsilon_{1}\epsilon_{2}}(\xi)\phi(\xi t,t)d\xi)dt
\cr\noalign {\vskip3truemm}
=\lim\limits_{\epsilon_{1},\epsilon_{2}\rightarrow0}\dis \int_{0}^{+\infty} t(\int_{-\infty}^{+\infty}H_{u}(\xi-\sigma)\phi(\xi t,t)d\xi)dt\cr\noalign {\vskip3truemm}
=\dis \int_{0}^{+\infty} \int_{-\infty}^{+\infty}H_{u}(\xi-\sigma)\phi(x,t)dxdt,
\end{array}
\end{align}
which means that
\begin{align}
\lim\limits_{\epsilon_{1},\epsilon_{2}\rightarrow0}\dis \int_{0}^{+\infty} \int_{-\infty}^{+\infty}(u^{\epsilon_{1}\epsilon_{2}}(x/t)
-H_{u}(\xi-\sigma))\phi(x,t)dxdt=0.
\end{align}
In a similar way, we have
\begin{align}
\begin{array}{l}
\lim\limits_{\epsilon_{1},\epsilon_{2}\rightarrow0}\dis \int_{0}^{+\infty} \int_{-\infty}^{+\infty}v^{\epsilon_{1}\epsilon_{2}}(x/t)\phi(x,t)
dxdt\cr\noalign{\vskip3truemm}
=\dis \int_{0}^{+\infty} \int_{-\infty}^{+\infty}H_{v}(\xi-\sigma)\phi(x,t)dxdt+\dis \int_{0}^{+\infty}(\sigma[v]-[uv])t\phi(\sigma t,t)dt,
\end{array}
\end{align}
that is
\begin{align}\label{eq5.24}
\lim\limits_{\epsilon_{1},\epsilon_{2}\rightarrow0}\dis \int_{0}^{+\infty} \int_{-\infty}^{+\infty}(v^{\epsilon_{1}\epsilon_{2}}(x/t)
-H_{v}(\xi-\sigma))\phi(x,t)dxdt=\dis \int_{0}^{+\infty}(\sigma[v]-[uv])t\phi(\sigma t,t)dt.
\end{align}

According to the definition \eqref{eq2.2}, the right side of \eqref{eq5.24} can be rewritten as
$$
\begin{array}{l}
\dis \int_{0}^{+\infty}(\sigma[v]-[uv])t\phi(\sigma t,t)dt
=\Big\langle w_{1}(\cdot)\delta_S,
\varphi(\cdot,\cdot)\Big\rangle,
\end{array}
$$
where
$$w_{1}(t)=t(\sigma[v]-[uv]).$$
The proof of Theorem 5.1 is completed. \ \ \ $\square$

\subsection{Formation of vacuum state for the system \eqref{eq1.1}}

Now we investigate the formation of vacuum state for the system \eqref{eq1.1}
when $(u_{+},v_{+})\in R_{1}R_{2}(u_{-},v_{-})$.
For this case, the solution of Riemann problem  \eqref{eq1.1} and \eqref{eq1.6} consists of
$R_{1}$ and $R_{2}$. Denoted the intermediate state
by $(u^{\epsilon_{1}\epsilon_{2}}_{\ast},v^{\epsilon_{1}\epsilon_{2}}_{\ast})$, then we have
\begin{align}\label{eq5.25}
R_{1}(u_{-},v_{-}):\left\{\begin{array}{l}
 \xi=\lambda^{\epsilon_{1}\epsilon_{2}}_{1} =u-\dfrac{1}{2}\epsilon_{2}-\dfrac{\sqrt{\epsilon^{2}_{2}+4\epsilon_{1}v^{2}}}{2}, \cr\noalign {\vskip2truemm}
u^{\epsilon_{1}\epsilon_{2}}_{\ast}-\dfrac{1}{2}\Big(-\sqrt{\epsilon^{2}_{2}+4\epsilon_{1}(v^{\epsilon_{1}\epsilon_{1}}_{\ast})^{2}}
+\epsilon_{2}\ln(\sqrt{\epsilon^{2}_{2}+4\epsilon_{1}(v^{\epsilon_{1}\epsilon_{2}}_{\ast})^{2}}+\epsilon_{2})\Big)=C_{1},
\cr\noalign {\vskip2truemm}
 u_{\ast}> u_{-}, \ \ v_{\ast}< v_{-}
\end{array}\right.
\end{align}
and
\begin{align}\label{eq5.26}
R_{2}(u_{-},v_{-}): \left\{\begin{array}{l}
 \xi=\lambda^{\epsilon_{1}\epsilon_{2}}_{2} =u-\dfrac{1}{2}\epsilon_{2}+\dfrac{\sqrt{\epsilon^{2}_{2}+4\epsilon_{1}v^{2}}}{2}, \cr\noalign {\vskip2truemm}
u^{\epsilon_{1}\epsilon_{2}}_{\ast}-\dfrac{1}{2}\Big(\sqrt{\epsilon^{2}_{2}+4\epsilon_{1}(v^{\epsilon_{1}\epsilon_{2}}_{\ast})^{2}}
+\epsilon_{2}\ln(\sqrt{\epsilon^{2}_{2}+4\epsilon_{1}(v^{\epsilon_{1}\epsilon_{2}}_{\ast})^{2}}-\epsilon_{2})\Big)=C_{2},\cr\noalign {\vskip2truemm}
 u_{\ast}< u_{+}, \ \ v_{\ast}<v_{+}.
\end{array}\right.
\end{align}

It is easy to find that
\begin{align}\label{eq5.27}
\begin{array}{l}
u_{+}-u_{-}=\dfrac{1}{2}\Big(-2\sqrt{\epsilon^{2}_{2}+4\epsilon_{1}(v^{\epsilon_{1}\epsilon_{2}})_{\ast}^{2}}+\sqrt{\epsilon^{2}_{2}+4\epsilon_{1}v_{-}^{2}}+\sqrt{\epsilon^{2}_{2}+4\epsilon_{1}v_{+}^{2}}
\cr\noalign{\vskip3truemm}
\hspace{2cm} +\epsilon_{2}\ln(\sqrt{\epsilon^{2}_{2}+4\epsilon_{1}(v^{\epsilon_{1}\epsilon_{2}}_{\ast})^{2}}+\epsilon_{2})-\epsilon_{2}\ln(\sqrt{\epsilon^{2}_{2}+4\epsilon_{1}(v^{\epsilon_{1}\epsilon_{2}}_{\ast})^{2}}-\epsilon_{2})
\cr\noalign{\vskip3truemm}
\hspace{2cm} -\epsilon_{2}\ln(\sqrt{\epsilon^{2}_{2}+4\epsilon_{1}v_{-}^{2}}+\epsilon_{2})+\epsilon_{2}\ln(\sqrt{\epsilon^{2}_{2}+4\epsilon_{1}v_{+}^{2}}-\epsilon_{2})\Big).
\end{array}
\end{align}
If $\lim\limits_{\epsilon_{1}, \epsilon_{2}\rightarrow0}v^{\epsilon_{1}\epsilon_{2}}_{\ast}=K\in(0,\mbox{min}(v_{-},v_{+}))$,
then we immediately get $u_{+}=u_{-}$ from \eqref{eq5.27}, which
contradicts with $u_{+}>u_{-}$. So we can conclude that  $\lim\limits_{\epsilon_{1}, \epsilon_{2}\rightarrow0}v_{\ast}=0$,
that is, the vacuum state will appear. In addition, from \eqref{eq5.25} and \eqref{eq5.26}, we have
 $\lim\limits_{\epsilon_{1}, \epsilon_{2}\rightarrow0}\lambda^{\epsilon_{1}, \epsilon_{2}}_{1}=\lim\limits_{\epsilon_{1}, \epsilon_{2}\rightarrow0}\lambda^{\epsilon_{1}, \epsilon_{2}}_{2}=u$.
At this moment, two rarefaction waves become two contact discontinuities $\xi=x/t=u_{\pm}$.
Therefore, we have the following theorem.

 \vspace{0.3cm}
\noindent{\bf Theorem 5.2.} Let $u_{+}>u_{-}$. For any fixed $\epsilon_{1}, \epsilon_{2}>0$, assume that $(u^{\epsilon_{1}\epsilon_{2}}, v^{\epsilon_{1}\epsilon_{2}})$ is a two-rarefaction wave
solution of the system \eqref{eq1.1} and \eqref{eq1.6}. Then, as
$\epsilon_{1}, \epsilon_{2}\rightarrow0$,
the two rarefaction waves become two contact discontinuities connecting the constant states $(u_{\pm},v_{\pm})$ and the vacuum
($v=0$), which form a vacuum solution of \eqref{eq1.2} and \eqref{eq1.6}.

\section{Limits of Riemann solutions to \eqref{eq1.1} when $\epsilon_{1}\rightarrow0$}

This section discusses the limit behaviors of Riemann solutions to \eqref{eq1.1} and \eqref{eq1.6} as $\epsilon_{1}\rightarrow0$.
Because our focus is the delta shock wave, so we first discuss the situation
$(u_{+},v_{+})\in III(u_{-},v_{-})$, that is, $u_{+}<u_{-}-2\epsilon_{2}$.

\vspace{0.3cm}
{\bf Lemma 6.1.} When $(u_{+},v_{+})\in III(u_{-},v_{-})$, there exists a positive parameter $\epsilon_{0}$ such that $(u_{+},v_{+})\in S_{1}S_{2}(u_{-},v_{-})$ when
$0<\epsilon_{1}<\epsilon_{0}$.

\vspace{0.3cm}
{\bf Proof.} \ \ When $(u_{+},v_{+})\in III(u_{-},v_{-})$, we have $u_{-}>u_{+}+2\epsilon_{2}$.
All states $(u,v)$ connected with $(u_{-},v_{-})$ by a backward shock wave $S_{1}$ or a forward shock wave
$S_{2}$ satisfy
\begin{align}
\begin{array}{l}
\dfrac{u-u_{-}}{v-v_{-}}=\dfrac{\epsilon_{2}-\sqrt{\epsilon^{2}_{2}+4\epsilon_{1}(v+v_{-})^{2}}}{v+v_{-}},\ \ u<u_{-},\ \ v>v_{-},
\cr\noalign {\vskip2truemm}
\dfrac{u-u_{-}}{v-v_{-}}=\dfrac{\epsilon_{2}+\sqrt{\epsilon^{2}_{2}+4\epsilon_{1}(v+v_{-})^{2}}}{v+v_{-}},\ \ u<u_{-},\ \ v<v_{-}.
\end{array}
\end{align}
If $v_{-}=v_{+}$, then we may take $\epsilon_{0}$ as any real positive number. Otherwise, by taking
$$(\dfrac{u_{+}-u_{-}}{v_{+}-v_{-}}(v_{-}+v_{+})-\epsilon_{2})^{2}=\Big(\sqrt{\epsilon^{2}_{2}+4\epsilon_{1}(v_{+}+v_{-})^{2}}\Big)^{2},$$
one can solve that
\begin{align}
\epsilon_{0}=\dfrac{\dfrac{u_{+}-u_{-}}{v_{+}-v_{-}}(v_{-}+v_{+})-\epsilon_{2})^{2}-\epsilon^{2}_{2}}{4(v_{-}+v_{+})^{2}}.
\end{align}
Considering that $u_{+}<u_{-}-2\epsilon_{2}$, so this lemma is ture. \ \ \ $\square$
\\

When $0<\epsilon_{1}<\epsilon_{0}$, besides two constant states $(u_{\pm},v_{\pm})$, the Riemann solution consists of a backward
shock wave $S_{1}$ and a forward shock wave $S_{2}$ with the intermediate state $(u^{\epsilon_{2}}_{\ast},v^{\epsilon_{2}}_{\ast})$
which is determined by
\begin{align}\label{eq6.3}
S_{1}(u_{-},v_{-}):\left\{\begin{array}{l}
 \sigma^{\epsilon_{2}}_{1} =u_{-}+\frac{v^{\epsilon_{2}}_{\ast}(u^{\epsilon_{2}}_{\ast}-u_{-})}{v^{\epsilon_{2}}_{\ast}-v_{-}}-\epsilon_{2}, \cr\noalign {\vskip2truemm}
u^{\epsilon_{2}}_{\ast}=u_{-}+(v^{\epsilon_{2}}_{\ast}-v_{-})\dfrac{\epsilon_{2}-\sqrt{\epsilon^{2}_{2}+4\epsilon_{1}(v^{\epsilon_{2}}_{\ast}+v_{-})^{2}}}{v^{\epsilon_{2}}_{\ast}+v_{-}},
\cr\noalign {\vskip2truemm}
 u^{\epsilon_{2}}_{\ast}<u_{-},
\ \ v^{\epsilon_{2}}_{\ast}>v_{-}
\end{array}\right.
\end{align}
and
\begin{align}\label{eq6.4}
S_{2}(u_{-},v_{-}):\left\{\begin{array}{l}
 \sigma^{\epsilon_{2}}_{2} =u_{+}+\frac{v^{\epsilon_{2}}_{\ast}(u_{+}-u^{\epsilon_{2}}_{\ast})}{v_{+}-v^{\epsilon_{2}}_{\ast}}-\epsilon_{2}, \cr\noalign {\vskip2truemm}
u_{+}=u^{\epsilon_{2}}_{\ast}+(v_{+}-v^{\epsilon_{2}}_{\ast})\dfrac{\epsilon_{2}+\sqrt{\epsilon^{2}_{2}+4\epsilon_{1}(v^{\epsilon_{2}}_{\ast}+v_{+})^{2}}}{v^{\epsilon_{2}}_{\ast}+v_{+}},
\cr\noalign {\vskip2truemm}
 u^{\epsilon_{2}}_{\ast}>u_{+},\ \ v^{\epsilon_{2}}_{\ast}>v_{+}.
\end{array}\right.
\end{align}

\vspace{0.2cm}
\noindent{ \bf Lemma 6.1} $\lim\limits_{\epsilon_{1}\rightarrow0}v^{\epsilon_{2}}_{\ast}=+\infty$.

\vspace{0.2cm}
\noindent{\bf Proof.} From the second equations of \eqref{eq6.3} and \eqref{eq6.4}, one can obtain
\begin{align}
u_{-}-u_{+}=(v_{\ast}-v_{+})\dfrac{\epsilon_{2}+\sqrt{\epsilon^{2}_{2}+4\epsilon_{1}(v_{\ast}+v_{+})^{2}}}{v_{\ast}+v_{+}}
+(v_{-}-v_{\ast})\dfrac{\epsilon_{2}-\sqrt{\epsilon^{2}_{2}+4\epsilon_{1}(v_{\ast}+v_{-})^{2}}}{v_{\ast}+v_{-}}.
\end{align}
If $\lim\limits_{\epsilon_{1}\rightarrow0}v_{\ast}=K\in(\mbox{max}({v_{-},v_{+}}),+\infty)$, then taking
$\epsilon_{1}\rightarrow0$ yields that
\begin{align}
u_{-}-u_{+}=2\epsilon_{2}\frac{K-v_{+}}{K+v_{+}}<2\epsilon_{2},
\end{align}
which contradicts with $u_{-}>u_{+}+2\epsilon_{2}$. Thus,
$\lim\limits_{\epsilon_{1}\rightarrow0}v^{\epsilon_{2}}_{\ast}=+\infty$. The proof is completed.  \ \ \ $\square$

Using the same way as used in section 5, we can easily obtain the following lemmas.

\vspace{0.3cm}
\noindent{ \bf Lemma 6.2} $\lim\limits_{\epsilon_{1}\rightarrow0}\sqrt{\epsilon^{2}_{2}+4\epsilon_{1}(v^{\epsilon_{2}}_{\ast}+v_{-})^{2}}=
\lim\limits_{\epsilon_{1}\rightarrow0}\sqrt{\epsilon^{2}_{2}+4\epsilon_{1}(v^{\epsilon_{2}}_{\ast}+v_{+})^{2}}=\dfrac{u_{-}-u_{+}}{2}$.

\vspace{0.3cm}
\noindent{ \bf Lemma 6.3} $ \lim\limits_{\epsilon_{1}\rightarrow0}u^{\epsilon_{2}}_{\ast}=\dfrac{u_{-}+u_{+}}{2}+\epsilon_{2},  \ \ \ \ \lim\limits_{\epsilon_{1}\rightarrow0}\sigma^{\epsilon_{2}}_{1}=\lim\limits_{\epsilon_{1}\rightarrow0}\sigma^{\epsilon_{2}}_{2}
=\dfrac{u_{-}+u_{+}}{2}.$

\vspace{0.3cm}
\noindent{\bf Lemma 6.4} $ \lim\limits_{\epsilon_{1}\rightarrow0}\dis\int_{\sigma^{\epsilon_{2}}_{1}}^{\sigma^{\epsilon_{2}}_{2}}v^{\epsilon_{2}}_{\ast}d\xi
=\frac{1}{2}(v_{+}(u_{-}-u_{+}+2\epsilon_{2})-v_{-}(u_{+}-u_{-}+2\epsilon_{2})).$

These lemmas show that, when $\epsilon_{1}\rightarrow0$, the velocities
of shock waves $S_{1}$ and $S_{2}$ approach to $\sigma^{\epsilon_{2}}$, which means that
$S_{1}$ and $S_{2}$ coincide. Besides, the intermediate state $v^{\epsilon_{2}}_{\ast}$ becomes
singular which determines the delta-shock solution of \eqref{eq1.5}  and \eqref{eq1.6}.
Similar to the proof in Theorem 5.1, we can draw the conclusion as follows.

\vspace{0.3cm}
\noindent{ \bf Theorem 6.1.} Let  $(u_{+},v_{+})\in III(u_{-},v_{-})$. For any fixed $\epsilon_{1}$, $\epsilon_{2}>0$, assume that $(u^{\epsilon_{1}\epsilon_{2}}, v^{\epsilon_{1}\epsilon_{2}})$ is a two-shock
Riemann solution of \eqref{eq1.1} and \eqref{eq1.6} constructed in Section 4. Then, as
$\epsilon_{1}\rightarrow0$, the limit of solution $(u^{\epsilon_{1}\epsilon_{2}}, v^{\epsilon_{1}\epsilon_{2}})$
is a delta shock wave connecting $(u_{\pm},v_{\pm})$, which is expressed by \eqref{eq3.8} with \eqref{eq3.10}  and \eqref{eq3.12}.

\vspace{0.3cm}
Next we discuss the limit behavior of the Riemann solution to \eqref{eq1.1}  and \eqref{eq1.6}
in the case $(u_{+},v_{+})\in I(u_{-},v_{-})$  as $\epsilon_{1}\rightarrow0$,
that is $u_{-}<u_{+}$.

\vspace{0.3cm}
{\bf Lemma 6.2.} When $(u_{+},v_{+})\in I(u_{-},v_{-})$, there exists a positive parameter $k_{0}$ such that $(u_{+},v_{+})\in R_{1}R_{2}(u_{-},v_{-})$ when $0<\epsilon_{1}<k_{0}$.

\vspace{0.3cm}
\noindent{\bf Proof.}  In this case, if $v_{+}<v_{-}$,
there exists a $k_{1}>0$, such that $(u_{+},v_{+})\in R_{1_{k_{1}}}$.
In fact, from  $(u_{+},v_{+})\in R_{1_{k_{1}}}$, we have
\begin{align}
u_{+}-u_{-}=\dis\int^{v_{+}}_{v_{-}}\dfrac{\epsilon_{2}-\sqrt{\epsilon^{2}_{2}+4\epsilon_{1}v^{2}}}{2v}\mbox{d}v, \ \ v_{+}<v_{-}.
\end{align}
Using the mean value theorem, we take
\begin{align}
k_{1}=\dfrac{(2v_{0}\frac{u_{+}-u_{-}}{v_{+}-v_{-}}-\epsilon_{2})^{2}-\epsilon^{2}_{2}}{4v_{0}^{2}}, \ \ v_{-}<v_{0}<v_{+}.
\end{align}
While if $v_{-}<v_{+}$,
there exists a $k_{2}>0$, such that $(u_{+},v_{+})\in R_{2_{k_{2}}}$.
In fact, from  $(u_{+},v_{+})\in  R_{2_{k_{2}}}$, we get
\begin{align}
u_{+}-u_{-}=\dis\int^{v_{+}}_{v_{-}}\dfrac{\epsilon_{2}+\sqrt{\epsilon^{2}_{2}+4\epsilon_{1}v^{2}}}{2v}\mbox{d}v, \ \ v_{-}<v_{+}.
\end{align}
By the mean value theorem, we take
\begin{align}
k_{2}=\dfrac{(2\bar{v_{0}}\frac{u_{+}-u_{-}}{v_{+}-v_{-}}-\epsilon_{2})^{2}-\epsilon^{2}_{2}}{4\bar{v_{0}}^{2}}, \ \ v_{-}<\bar{v_{0}}<v_{+}.
\end{align}
Taking $k_{0}=\min\{k_{1},k_{2}\}$, we get the result.  \ \ \ $\square$

When $0<\epsilon_{1}<k_{0}$, the Riemann solution to \eqref{eq1.1}  and \eqref{eq1.6}
 just consists of two rarefaction waves.
Letting $\epsilon_{1}\rightarrow0$ in \eqref{eq4.4} and \eqref{eq4.5}, $R_{1}$ and $R_{2}$
become the contact discontinuity $J$ and rarefaction wave $R$, which can be
expressed as
\begin{align}
J:\left\{\begin{array}{l}
 \xi =\lambda^{\epsilon_{2}}_{1}=u-\epsilon_{2}, \cr\noalign {\vskip2truemm}
 \mbox{d}u=0, \ \ \ u^{\epsilon_{2}}_{\ast}=u_{-}
\end{array}\right.
\end{align}
and rarefaction wave
\begin{align}
R(u_{-},v_{-}):\left\{\begin{array}{l}
 \xi =\lambda^{\epsilon_{2}}_{2}=u, \cr\noalign {\vskip2truemm}
 v\mbox{d}u=\epsilon_{2}\mbox{d}v, \ \ \ u^{\epsilon_{2}}_{\ast}<u_{+},
\end{array}\right.
\end{align}
where $(u^{\epsilon_{2}}_{\ast},v^{\epsilon_{2}}_{\ast})$ is the intermediate state
of $J$ and $R$. After a simple calculation, we can get
\begin{align}
(u^{\epsilon_{2}}_{\ast},\ v^{\epsilon_{2}}_{\ast})=(u_{-},\ v_{+}\exp(\frac{u_{-}-u_{+}}{\epsilon_{2}})).
\end{align}
Then, we have the following theorem.

\vspace{0.3cm}
\noindent{ \bf Theorem 6.2.} Let  $(u_{+},v_{+})\in I(u_{-},v_{-})$. For any fixed $\epsilon_{1}$, $\epsilon_{2}>0$, assume that $(u^{\epsilon_{1}\epsilon_{2}}, v^{\epsilon_{1}\epsilon_{2}})$ is a two-rarefaction wave
Riemann solution of the system \eqref{eq1.1} and \eqref{eq1.6} constructed in Section 4.
Then as $\epsilon_{1}\rightarrow0$,  the limit of solution $(u^{\epsilon_{1}\epsilon_{2}}, v^{\epsilon_{1}\epsilon_{2}})$
is a contact discontinuity and a rarefaction wave connecting $(u_{\pm},v_{\pm})$.

\vspace{0.3cm}
We have proven that when  $u_{-}>u_{+}+2\epsilon_{2}$ or $u_{-}<u_{+}$,
the solutions to the Riemann problem of \eqref{eq1.1} and \eqref{eq1.6} just are the solutions to the
Riemann problem for \eqref{eq1.5} with the same initial data as $\epsilon_{1}\rightarrow0$.
The same conclusions are true for the rest cases, and we omit the discussions.

\noindent{\bf Remark}. The processes of formation of delta shock waves and vacuum states can be examined
with some numerical results as $\epsilon_{1}$ and $\epsilon_{2}$ decrease. The numerical simulations will be presented
in the version for publication.

\end{document}